\documentclass[a4paper,11pt]{article}
\usepackage[english]{babel}
\usepackage{theorem}
\usepackage{amssymb}
\usepackage{amsfonts}
\usepackage{amsmath}
\usepackage[all]{xy}
\usepackage{geometry}
\usepackage{ae,aecompl}
\usepackage{eurosym}
\usepackage{verbatim}

\newtheorem{thm}{Theorem}[section]
\newtheorem{prop}[thm]{Proposition}
\newtheorem{lem}[thm]{Lemma}
\newtheorem{cor}[thm]{Corollary}

\newtheorem{remk}[thm]{Remark}

\makeatletter
\begingroup
\gdef\th@upshape{\normalfont
  \def\@begintheorem##1##2{%
        \item[\hskip\labelsep \theorem@headerfont ##1\ ##2.]}%
\def\@opargbegintheorem##1##2##3{%
   \item[\hskip\labelsep \theorem@headerfont ##1\ ##2\ (##3).]}}
\endgroup
\makeatother

\theoremstyle{upshape}

\newtheorem{rem}[thm]{Remark}
\newtheorem{defn}[thm]{Definition}

\title{{\bf Noncommutative Complex Structures on Quantum Homogeneous Spaces}}
\author{R\'{e}amonn \'{O} Buachalla \footnote{Supported by the Grant GACR P201/12/G028}}
\date{}

\date{}

\def\bal#1\eal{\begin{align}#1\end{align}}
\def\bas#1\eas{\begin{align*}#1\end{align*}}
\def\bit{\begin{itemize}}
\def\eit{\end{itemize}}
\def\ed{\end{document}}
\def\bet{\begin{enumerate}}
\def\eet{\end{enumerate}}

\def\a{\alpha}

\def\g{\gamma}
\def\d{\delta}
\def\e{\varepsilon}
\def\f{\varphi}

\def\l{\lambda}

\def\r{\rho}
\def\s{\sigma}

\def\t{\tau}
\def\w{\omega}

\def\Om{\Omega}

\def\del{\partial}
\def\adel{\ol{\partial}}
\def\DEL{\Delta}
\def\G{\Gamma}

\def\P{\Phi}

\def\bC{{\mathbb C}}
\def\bN{{\mathbb N}}

\def\bZ{{\mathbb Z}}

\def\E{{\cal E}}
\def\F{{\cal F}}

\def\T{{\cal T}}

\def\exd{\mathrm{d}}
\def\demo{\noindent \emph{\textbf{Proof}.\ }~}

\def\unit{\mathrm{U}}
\def\counit{\mathrm{C}}
\def\frame{\mathrm{U}}
\def\frameadj{\mathrm{C}}

\def\id{\mathrm{id}}
\def\ker{\mathrm{ker}}
\def\proj{\mathrm{proj}}

\def\spn{\mathrm{span}}
\def\th{^\mathrm{th}}

\def\setl{ ~ \text{\Big{|}} ~ }
\def\hol{^{(1,0)}}
\def\ahol{^{(0,1)}}

\def\inv{^{-1}}

\def\coby{\, \square_{H}}
\def\oby{\otimes}

\def\wed{\wedge}
\def\sseq{\subseteq}
\def\tl{\triangleleft}

\def\wt{\widetilde}

\def\ol{\overline}
\def\la{\left\langle}
\def\\la{\left\langle}
\def\ra{\right\rangle}
\def\>{\right\rangle}
\def\bs{\backslash}
\def\mto{\mapsto}

\def\qed{\hfill\ensuremath{\square}\par}
\def\qf3{\bC_q[F_3]}

\def\csu2{\bC_q[SU_2]}
\def\cs2{\bC_q[S^2]}
\def\csun{\bC_q[SU_N]}

\def\ccp1{\bC P^{1}}
\def\cp1{\bC_q[\bC P^{1}]}
\def\cpn{\bC_q[\bC P^{N-1}]}

\def\ccpn{\bC P^{N-1}}

\def\qf3{\bC_q[F_3]}
\def\wqf3{\Om^1_q[F_3]}
\def\qsu3{\bC_q[SU_3]}
\def\wqsu3{\Om^1_q[SU_3]}

\def\usl2{\mathcal{U}(\mathfrak{sl}(2))}
\def\wsun{\Om^1_q(SU_N)}
\def\ws2{\Om^1_q(S^2)}
\def\wsu2{\Om^1_q(SU_2)}
\def\wcpn{\Om^1_q(\ccpn)}

\def\hol{^{(1,0)}}
\def\ahol{^{(0,1)}}

\def\n2{_{(-2)}}
\def\m2{_{(-2)}}
\def\m1{_{(-1)}}
\def\0{_{(0)}}
\def\1{_{(1)}}
\def\2{_{(2)}}
\def\3{_{(3)}}
\def\4{_{(4)}}
\def\5{_{(5)}}
\def\hol{^{(1,0)}}
\def\ahol{^{(0,1)}}

\def\cvect{{}_{\bC} \hspace{-.030cm}\mathrm{Mod}}

\def\mm{{}_M \hspace{-.030cm}\mathrm{Mod}}

\def\gmm{{}^{G}_M \hspace{-.030cm}\mathrm{Mod}}
\def\gmmm{{}^{G}_M \hspace{-.030cm}\mathrm{Mod}_M}
\def\mh{\mathrm{Mod}^H}
\def\mhm{\mathrm{Mod}^H_M}
\def\lgmmm{{}^{G}_M \hspace{-.030cm}{\mathrm{Mod}}_0}
\def\lmhm{{\mathrm{Mod}}^H_0}
\def\ggmg{{}^{G}_G \hspace{-.030cm}{\mathrm{Mod}}_G}

\def\acs{almost complex structure~}
\def\alg{algebra~}
\def\algn{algebra}
\def\algs{algebras~}

\def\ff{faithfully flat~}

\def\fodc{first-order differential calculus~}
\def\fodcn{first-order differential calculus}

\def\hk{Heckenberger--Kolb~}

\def\hakn{Heckenberger and Kolb}

\def\lc{left-covariant~}

\def\nccg{noncommutative complex geometry~}
\def\nccgn{noncommutative complex geometry}

\def\ncgn{noncommutative geometry}
\def\nc{noncommutative~}

\def\qhs{quantum homogeneous space~}
\def\qhsn{quantum homogeneous space}

\def\iff{if and only if~}

\def\st{such that~}
\def\stn{such that}

\def\wrt{with respect to~}

\DeclareMathOperator{\dt}{det}

\def\mpr{maximal prolongation~}
\def\mprn{maximal prolongation}

\def\uqsl2{U_q(\frak{sl_2)}}
\def\tuqsl2{\wt{U}_q(\frak{sl}_2)}
\def\tu1sl2{\wt{U}_1(\frak{sl}_2)}

\begin{document}

\maketitle

\begin{abstract}
A new framework for \nc complex geometry on quantum homogeneous spaces is introduced. The main ingredients used are covariant differential calculi and Takeuchi's categorical equivalence for quantum homogeneous spaces. A number of basic results are established, producing a simple set of necessary and sufficient conditions for  \nc complex structures to exist. Throughout, the framework is applied to the quantum projective spaces endowed with the Heckenberger--Kolb calculus.
\end{abstract}

\section{Introduction}

Classical complex geometry is a subject of remarkable richness and beauty with deep connections to modern  physics. Yet despite over twenty-five years of \ncgn, the development of \nc  complex geometry is still in its infancy. What we do have is a large number of examples which demand consideration as \nc complex spaces. These include, amongst others, noncommutative tori \cite{DiengSchw}, noncommutative projective algebraic varieties \cite{SV}, fuzzy flag manifolds \cite{DIAS}, and (most importantly from the point of view of this paper) examples arising from the theory of quantum groups \cite{HKdR,Maj}. 

Thus far, there have been two attempts to formulate a general framework for \nc complex geometry. The first, due to Khalkhali, Landi, and van Suijlekom \cite{KLVSCP1}, was introduced to provide a context for their work on the \nccg of the Podle\'s sphere. This  followed on from earlier work of Majid \cite{Maj},   Schwartz and Polishchuk \cite{PolishSchw}, and  Connes \cite{ConnesCuntz, CONN}. Khalkhali and Moatadelro  \cite{KKCP2,KKCPN} would  go on to apply this framework to D'Andrea and D\c{a}browski's work \cite{DDCPN} on the  higher order quantum projective spaces.
 
Subsequently, Beggs and Smith introduced a second  more comprehensive approach to \nccg  in  \cite{EBPS}. Their motive was to provide a framework for quantising the intimate relationship between complex differential geometry and complex projective geometry. They foresee that the rich interaction between algebraic and analytic techniques occurring in the classical setting will carry over to the \nc world.

The more modest aim of this paper is to begin the development of a theory of \nccg for quantum group homogeneous spaces. This will be done very much in the style of Majid's noncommutative Riemannian geometry \cite{Maj,MajPrimer}, with the only significant difference being that here we will not need to assume that our quantum homogeneous spaces are Hopf--Galois extensions, while we will assume a faithful flatness property. This assumption allows us to use Takeuchi's categorical equivalence to establish a simple set of necessary and sufficient conditions for covariant  complex structures to exist. In subsequent work, it is intended to build upon these results and formulate \nc generalisations of Hodge theory and K\"ahler geometry for quantum homogeneous spaces \cite{MMF3}. Indeed, the first steps in this direction have already been taken \cite{MMFPhD}.
  
For this undertaking to be worthwhile, however,  it will need to be applicable to a good many interesting examples. Recall that classically one of the most important classes of homogeneous complex manifolds is the family of generalised flag manifolds. As has been known for a long time, these spaces admit a direct $q$-deformation in terms of the Drinfeld--Jimbo quantum groups \cite{LR91,Soibel2,TT}. Somewhat more recently, it was shown by Heckenberger and Kolb \cite{HKdR} that the Dolbeault double complex of the irreducible flag manifolds survives this $q$-deformation intact. This result gives us one of the most important families of  \nc complex structures that we have, and as such, provides an invaluable testing ground for any newly proposed theory of  \nccgn.


In this paper we show that, for the special case of  quantum projective space, the work of Heckenberger and Kolb can be understood in terms of our general framework for \nccgn. This allows for a significant simplification of the required calculations, and helps identify some of the underlying general processes at work. It is foreseen that this work will prove easily extendable to all the irreducible quantum flag manifolds. Moreover, it is hoped to extend it even further to include all the quantum flag manifolds,  and in so doing, produce  new examples of \nc complex structures. 

\bigskip

The paper is organised as follows: 
In section 2 some well-known material about quantum homogeneous spaces, Takeuchi's categorical equivalence, covariant differential calculi, almost complex structures, and  complex structures is recalled.

In Section 3  the quantum special unitary group, and the quantum projective spaces, as well as the Heckenberger--Kolb calculus for these spaces, is discussed. 

In Section 4 one of the basic results of the paper {\bf Proposition \ref{monoidal.equiv}} is established: It shows that for a special subcategory of $\mhm$, the monoidal structure induced on it by the canonical monoidal structure of $\gmmm$ (through Takeuchi's equivalence)  is equivalent to the vector space tensor product.

In Section 5, {\bf Theorem \ref{theframingtheorem}}  shows how to find an explicit description of the \mpr 
of a covariant \fodc in terms of a certain  ideal $I^{(1)} \sseq M^+$.

In Section 6 the notion of factorisability for almost complex structures is introduced, and  a simple set of necessary and sufficient conditions for factorisable almost complex structures to exist is established.

Finally, in Section 7, {\bf Proposition \ref{prop-integrability}} gives a simple method  for verifying that an almost complex structure is a complex structure.


Throughout, the family of quantum projective spaces, endowed with the Heckenberger--Kolb calculus, is taken as the motivating set of examples. In each section, the newly constructed general theory is applied to these examples in detail, building up to an explicit presentation of their $q$-deformed Dolbeault double complexes.

\section{Preliminaries and First Results}

In this section we recall Takeuchi's categorical equivalence for quantum homogeneous spaces, some of its applications to the theory of covariant differential calculi, and finally the definition of a complex structure. 

 \subsection{Quantum Homogeneous Spaces}

Let $G$ be a Hopf  \alg with comultiplication $\DEL$, counit $\e$, antipode $S$, unit $1$, and multiplication $m$. 
Throughout, we use Sweedler notation, as well as denoting $g^+ := g - \e(g)1$, for $g \in G$, and $V^+ = V \cap \ker(\e)$, for $V$   a subspace of $G$. For a right $G$-comodule $V$ with coaction $\DEL_R$, we say that an element $v \in V$ is {\em coinvariant} if $\DEL_R(v) = v \oby 1$. We denote the subspace of all coinvariant elements by $V^G$, and call it the {\em coinvariant subspace} of the coaction. 
For $H$ a Hopf \algn, a {\em homogeneous} right $H$-coaction on $G$ is a coaction of the form $(\id \oby \pi) \circ \DEL$, where $\pi: G \to H$ is a Hopf \alg map. The coinvariant subspace of such a coaction is a sub\alg \cite[Proposition 1]{Tak}. 

\begin{defn}
We call the coinvariant sub\alg $M:=G^H$ of  a homogeneous coaction a {\em quantum homogeneous space} if $G$ is faithfully flat as a right module over $M$, which is to say  if  the functor  $G \oby_M -:\mm \to \cvect$, from the category of left $M$-modules to the category of complex vector spaces, maps a sequence to an exact sequence if and only if the original sequence is exact.
\end{defn}

In this paper we will {\em always} use the symbols $G,H,\pi$ and $M$ in this sense. We  also note that $G$ is itself a trivial example of a \qhsn, where $\pi = \e$. Moreover, the coproduct of $G$ restricts to a right $G$-coaction on $M$, and
\bal \label{piofm0}
\pi(m) = \e(m)1_H,  & & \text{for all } m \in M.
\eal
If $G$ and $H$ are Hopf $*$-\algn s, and $\pi$ is a Hopf $*$-\alg map, then $M$ is a $*$-sub\alg of $G$.


\subsection{Some Categories}\label{Takeq}

We now define the abelian categories $\gmmm$ and $\mhm$. The objects in $\gmmm$ are $M$-bimodules $\E$ (with left and right actions denoted by juxtaposition) endowed with a left $G$-coaction $\DEL_L$ \st 
\bal \label{gmmmdefn}
\DEL_L(m e m') =  m\1 e\m1 m'\1 \oby m\2 e\0 m'\2, & & \text{for all ~} m,m' \in M, e \in \E.
\eal
The morphisms in $\gmmm$ are the $M$-bimodule homomorphisms  that are also homomorphisms of left $G$-comodules. The objects in $\mhm$ are right $M$-modules $V$ (with right action denoted by $\tl$) endowed with a right $H$-coaction $\DEL_R$ such that 
\bal \label{mhm-compt}
\DEL_R(v \tl  m) =  v\0 \tl m\2 \oby S(\pi(m\1))v\1, & &  \text{for all ~} v \in V, m \in M.
\eal
The morphisms in $\mhm$ are the $M$-module homomorphisms that are also homomorphisms of right $H$-comodules.


\bigskip

Next we introduce a subcategory of $\gmmm$, and a subcategory of $\mhm$, that play important roles in the paper. The definition of the latter requires the following technical lemma.
\begin{lem}  For $\mh$ the category of right $H$-comodules, we have a fully faithful embedding 
\bal \label{embed} 
\mh \to \mhm, & & V \mto (V, \tl),
\eal
where $\tl$ is the the trivial right $M$-module structure, $v \tl  m = \e(m)v$, for $v \in V, m \in M$. 
\end{lem}
\demo
To show that $(V,\tl)$ is well-defined as an object in $\mhm$, we need to show that (\ref{mhm-compt}) is satisfied.  This is implied by (\ref{piofm0}) as follows:
\bas
v\0 \tl m\2 \oby S(\pi(m\1))v\1 & = v\0\e(m\2) \oby \e(m\1)v\1 = \e(m) v\0 \oby v\1 \\
                                               & = \DEL_R(\e(m) v) = \DEL_R(v \tl m).
\eas
Moreover, since any comodule map between $V$ and $W$ is trivially a module map \wrt $\tl$, it is clear that (\ref{embed}) defines a fully faithful functor. \qed

\begin{defn}
Denote by $\lgmmm$ the full subcategory of $\gmmm$ whose objects $\E$  satisfy $\E M^+ \sseq M^+\E$, and denote by $\lmhm$ the image of $\mh$ under the embedding in (\ref{embed}). 
\end{defn}

\subsection{Takeuchi's Categorical Equivalence}

If $\E \in \gmmm$, then $\E/(M^+\E)$ becomes an object in $\mhm$ with the obvious right $M$ action, and the right $H$-coaction 
\bal \label{comodstruc0}
\DEL_R(\ol{e}) = \ol{e\0} \oby S(\pi(e\m1)), & & e \in \E.
\eal
We define a functor $\Phi:\gmmm \to \mhm$  as follows:  
$\Phi(\E) :=  \E/(M^+\E)$, and if $g : \E \to \F$ is a morphism in $\gmmm$, then $\Phi(g):\Phi(\E) \to \Phi(\F)$ is the map to which $g$ descends on $\Phi(\E)$.

If $V \in \mhm$, then $G \coby V := (G \oby V)^H$  (where $G\oby V$ has the standard tensor product $H$-comodule structure  \cite[\textsection 1.3.2]{KSLeabh}) becomes an object in $\gmmm$ with $M$-bimodule structure
\begin{align*}
m \Big(\sum_i g^i \oby v^i\Big)  = \sum_i m g^i \oby v^i, & & \Big(\sum_i g^i \oby v^i\Big) m = \sum_i g^i m\1 \oby (v^i \tl  m\2),
\end{align*}
and left-$G$-coaction 
\begin{align*}
\DEL_L\Big(\sum_i g^i \oby v^i\Big) = \sum_i g^i\1 \oby g^i\2 \oby v^i.
\end{align*}
We define a functor $\Psi:\mhm \to \gmmm$ as follows:
$
\Psi(V) := G \coby V,
$
and if $\g$ is a morphism in $\mhm$, then $\Psi(\g) := \id \oby \g$.

\begin{thm} \cite[Theorem 1]{Tak} An equivalence of categories
between $\gmmm$ and $\mhm$, which we call {\em Takeuchi's equivalence},  is given by the functors $\Phi$ and $\Psi$ and the natural transformations
\begin{align}
\counit:\Phi \circ \Psi(V) \to V, & & \ol{\sum_i g^i \oby v^i} \mto \sum_{i} \e(g^i)v^i \label{counit},\\
\unit: \E \to \Psi \circ \Phi(\E), & & e \mto e\m1 \oby \ol{e\0}. \label{unit}
\end{align}
\end{thm}
We define the {\em dimension} of an object $\E \in \gmmm$ to be the vector space dimension of $\Phi(\E)$.

We now present an explicit formula for the inverse of $\unit$. We do so in a number of steps, so as to highlight some results that will be of use to us later.

\begin{lem}\cite[\textsection 1]{Tak}
An isomorphism is given by 
\bas
u: G\oby_M \E \to G \oby \Phi(\E), & & g \oby_M e \mto ge\m1 \oby \ol{e\0}.
\eas 
Moreover, the inverse of $u$ acts according to  $u\inv(g \oby \ol{e}) \mto  g S(e \m1) \oby_M e\0$.
\end{lem}


\begin{cor}
For $\E \in \gmmm$, the following diagram is commutative
\bas
\xymatrix{ 
\E   \ar[d]_{\unit}                 \ar[rrr]^{1 \oby_M \id}           & & &   G \oby_M \E \ar[d]^{u} \\
\Psi \circ \Phi(\E)   \,\,   \ar@{^{(}->}[rrr]                     & & & G \oby \Phi(\E),\\
}
\eas
where the inclusion in the bottom row is the obvious one. Hence, $1 \oby_M \id$ is an embedding. 
\end{cor}
\demo

It follows directly from the definitions of $\unit$ and $u$ that the diagram is commutative. Hence, since $u$ is an isomorphism, $1 \oby_M \id$  must be an embedding.
\qed

\begin{cor} \label{unitinv}
The inverse of $\unit$ is  given by
\bas 
\unit \, \inv: \Psi \circ \Phi(\E) \to G \oby_M \E,  &  & \sum_i g^i \oby \ol{e^i} \mto \sum_i  g^i S(e^i\m1) \oby_M e^i\0.
\eas
\end{cor}

\bigskip

Finally, we turn to the question of how Takeuchi's equivalence behaves upon restricting to the two subcategories $\lgmmm$ and $\lmhm$.


\begin{lem}
Takeuchi's equivalence restricts to an equivalence between the subcategories $\lgmmm$ and $\lmhm$.
\end{lem}
\demo
If $\E$ is an object in $\lgmmm$, then for any $e \in \E$, and $m \in M^+$, the fact that $\E M^+ \sseq M^+\E$ implies that $\ol{e} \tl m =0$. Hence, for any $n \in M$, we  have
\bas
\ol{e}\tl n & = \ol{e}\tl \left(n^+ + \e(n)1\right)  = \ol{e}\tl n^+ + \ol{e}\tl\left(\e(n)1\right) = \e(n)\ol{e},
\eas 
showing us that $\Phi(\E)$ is well-defined as an object in $\lmhm$. Conversely, if $V$ is an object in $\lmhm$, then for any element $\sum_i f^i \oby v^i$ in $\Psi(V)$,
\bas
\Big{(}\sum_i f^i \oby v^i\Big{)}m = \sum_i f^im\1 \oby \left(v^i \tl m\2\right) =  \sum_i f^i  m\1 \oby \e(m\2) v^i  = \sum_i f^im \oby v^i.
\eas
If $m \in M^+$, then $\sum_i f^im \oby v^i \in \ker(\frameadj)$.  But  $\ker(\frameadj) = M^+\Psi(V)$ so  $\left(\sum_i f^i \oby v^i\right)m \in M^+\Psi(V)$. Hence $\Psi(V)$ belongs to $\lgmmm$. This establishes the second assertion of the lemma.
\qed

\begin{remk} \label{rem2}
Roughly speaking, we view $\gmmm$ as generalising the category of equivariant vector bundles over a homogeneous space; $\mhm$ as generalising the category of representations of the isotropy subgroup; and Takeuchi's adjunction as generalising the well known equivalence between these categories \cite[\textsection 1]{Segal}.
\end{remk} 

\subsection{First-Order Differential Calculi}

Let $A$ be a unital algebra (in what follows all \algs are assumed to be unital). A {\em first-order differential calculus} over $A$ is a pair $(\Om^1,\exd)$, where $\Omega^1$ is an $A$-$A$-bimodule and $\exd: A \to \Omega^1$ is a linear map for which the {\em Leibniz rule} holds
\begin{align*}
\exd(ab)=a(\exd b)+(\exd a)b,&  & a,b,\in A,
\end{align*}
and for which $\Om^1 = \spn_{\bC}\{a\exd b\, | \, a,b \in A\}$. (Where no confusion arises, we will drop explicit reference to $\exd$ and denote a calculus by its bimodule $\Om^1$ alone.) We call an element of $\Om^1$ a {\em one-form}. An {\em isomorphism} between two first-order differential calculi $(\Om^1(A),\exd_{\Om})$ and $(\G^1(A),\exd_{\G})$ is a bimodule isomorphism $\f: \Om^1(A)  \to \Gamma^1(A)$ \st $\f \circ \exd_{\Om} = \exd_{\Gamma}$. The {\em direct sum} of two first-order differential calculi $(\Om^1(A),\exd_{\Om})$ and $(\G^1(A),\exd_{\G})$ is the calculus $(\Om^1(A) \oplus \Gamma^1(A), \exd_\Om + \exd_\G)$.

The {\em universal first-order differential calculus} over $A$ is the pair
$(\Om^1_u(A), \exd_u)$, where $\Om^1_u(A)$ is the kernel of the product map $m: A \otimes A \to A$ endowed
with the obvious bimodule structure, and $\exd_u$ is defined by
\begin{align*}
\exd_u: A \to \Omega^1_u(A), & & a \mto 1 \oby a - a \oby 1.
\end{align*}
By \cite[Proposition 1.1]{Wor}, every first-order differential calculus over $A$ is of the form $\left(\Omega^1_u(A)/N, \,\proj \circ \exd_u\right)$, where $N$ is a $A$-sub-bimodule of $\Omega^1_u(A)$, and  $\proj:\Omega^1_u(A) \to \Omega^1_u(A)/N$ is the canonical projection. Moreover, this gives a bijective correspondence between calculi and sub-bimodules. 

We say that a differential calculus $\Omega^1(M)$, over a \qhs $M$, is {\em covariant} if there exists a (necessarily unique) map $\DEL_L: \Om^1(M) \to G \oby \Om^1(M)$,  \st 
\bas
\DEL_L(m\exd n) = \DEL(m) (\id \oby \exd)\DEL(n),  &  & m,n \in M.
\eas 
Any covariant calculus $\Om^1(M)$ is naturally an object in $\gmmm$. Moreover, the universal calculus over $M$ is covariant, and covariance of any $\Om^1(M) \simeq \Om^1_u(M)/N$ is equivalent to $N$ being a sub-object of $\Om^1_u(M)$ in $\gmmm$. (Note that $\exd$ is not a morphism in $\gmmm$.)

The following theorem is a special case of more general results originally established by Hermisson \cite[Theorem 2]{Herm}, and Majid \cite[Theorem 2.1]{Maj}. 

\begin{thm} \label{MHW}
For a  \qhs $M$, considering $M^+$ as an object in $\mhm$ according to its obvious right $M$-module structure, and the right $H$-comodule structure $\DEL_{R}(m) = m\2 \oby S(\pi(m\1))$, for $m \in M^+$, it holds that:
\begin{enumerate}
\item
Covariant first-order differential calculi over $M$ are in bijective correspondence with sub-objects of $M^+$.

\item The sub-object corresponding to the calculus $\Om^1(M)$ is 
\bas
I^{(1)} := \Big\{  \sum_i \e(m_i)m^+_i \,  ~ \text{\Big{|}} ~  \, \sum_i m_i \exd n_i = 0\ \Big{\}}.
\eas
\item
Denoting $V^1:= M^+/I^{(1)}$,  (which we call the {\em cotangent space} of $\Om^1(M)$) we have an isomorphism
\bas
\s:\Phi\left(\Om^1(M)\right) \to V^1, & & \ol{ m \exd n} \mto \ol{\e(m)m^+}.
\eas
\end{enumerate}
\end{thm}
\demo
Applying the functor $\Phi$ to the collection of  sub-objects of $\Om^1_u(M)$ gives a correspondence between \lc calculi over $M$ and sub-objects of  $\Phi(\Om^1_u(M))$. The theorem then follows from the easily verifiable fact that an isomorphism is given by
\bas
\Phi\left(\Om^1_u(M)\right) \to M^+, & & \ol{m \exd n} \mto \e(m)n^+.
\eas
\qed

For the special case of a trivial \qhsn, this result reduces to  Woronowicz's  celebrated theorem classifying \lc calculi over a Hopf \alg $G$ \cite[Theorem 1.5]{Wor}. For such a calculus $\Om^1(G)$, we follow the standard convention of denoting its cotangent space by $\Lambda$, and calling it  the {\em space of left-invariant $1$-forms} of the calculus. 

\bigskip

We note that for $\Om^1(G)$ any calculus on $G$, the bimodule 
 $\Om^1(M) := \left\{ m \exd n \,|\, m,n \linebreak \in M\right\}$
has the natural structure of a \fodc over $M$. We call it the {\em restriction} of $\Om^1(G)$  to $M$.

\subsection{Differential Calculi}

For $(S,+)$ a commutative semigroup, an  {\em {$S$-graded \alg}} is an \alg $A$ equipped with a decomposition $A = \bigoplus_{s \in S}A^s$, where each $A^s$ is a linear
subspace of $A$, and $A^sA^t \sseq A^{s+t}$, for all $s,t
\in S$. If $a \in A^s$, then we say that $a$ is a {\em homogeneous element of degree $s$}.  A {\em homogeneous mapping of degree $t$} on $A$ is a linear mapping $L:A \to A$ \st if $a \in
A^s$, then $L(a) \in A^{s+t}$. We say that a subspace $B$ of $A$ is {\em homogeneous} if it admits a decomposition  $B = \bigoplus_{s \in S} B^s$, with $B^s \sseq A^s$, for all $s \in S$.

A pair $(A,\exd)$ is called a {\em complex} if $A$ is an ${\bN}_0$-graded \algn, and $\exd$ is a homogeneous mapping  of degree $1$ \st $\exd^2 = 0$. A triple $(A,\del,\ol{\del})$ is called a {\em double complex} if $A$ is an $\bN^2_0$-graded \algn, $\del$ is homogeneous mapping of degree $(1,0)$, $\ol{\del}$ is homogeneous mapping of degree $(0,1)$, and 
\bas
\del^2 = \ol{\del}^2 = 0, & & \del \circ \ol{\del} =  - \ol{\del} \circ \del.
\eas 
Observe that we can associate to any double complex $(A,\del,\ol{\del})$ the complex $(A,\del + \ol{\del})$.

A complex $(A,\exd)$ is called a {\em differential graded \alg} if $\exd$ is a {\em graded derivation}, which is to say, if it satisfies the {\em graded Leibniz rule}
\bas
\exd (ab) = \exd (a)b+(-1)^n a \exd b,       &  & \textrm{  for all $a \in A^n$, $b \in A$}.
\eas
The operator $\exd$ is called
the {\em differential} of the differential graded \algn. 
\begin{defn}
A {\em differential calculus} over an \alg $A$ is a differential
\alg  $(\Om(A),\exd)$ \st $\Om^0=A$, and
\begin{align} \label{pofu}
\Om^k = \spn_{\bC}\{a_0\exd a_1\wed  \cdots \wed \exd a_k \,|\, a_0, \ldots, a_k \in A\}.
\end{align}
\end{defn}
We use $\wedge$ to denote the multiplication between elements of a differential calculus when both are of order greater than or equal to $1$, otherwise we use juxtaposition. 

For any $A$-$A$-bimodule $\E$, we denote 
$
\T(\E) := \bigoplus_{k=0}^\infty \E^{\oby_A k}.
$
Endowed with the obvious structure of a graded \algn, we call $\T(\E)$ the {\em tensor \alg} of $\E$. 
Any first-order differential calculus $\Om^1(A)$ can be extended to a  differential calculus: For $N \sseq \Om^1_u(A)$ the sub-bimodule corresponding to $\Om^1(A)$, denote  
\bal \label{abused}
\Om^{\bullet}(A):= \T(\Om^1(A))/\la \exd N \ra,
\eal 
where $\la \exd N \ra$ is the ideal of  $\T(\Om^1(A))$ generated by $\exd N$, and by abuse of notation, $\exd N$ is the image in $(\Om^1(A))^{\oby_A 2})$ of $d_u \oby \exd_u(N)$ under the canonical projection $(\Om^1_u(A))^{\oby_A 2} \to (\Om^1(A))^{\oby_A 2}$.  The exterior derivative $\exd$ is easily seen to have a unique extension $\exd:\Om^{\bullet}(A) \to \Om^{\bullet}(A)$ \st  $(\Om^{\bullet}(A),\exd)$ is a  differential calculus. We call this  differential calculus the {\em \mpr} of $(\Om^1(A),\exd)$. The \mpr is unique in the sense that any other calculus extending $(\Om^1(A),\exd)$ can be obtained as a quotient of the \mprn.

If $\Om^1(M)$ is a covariant first-order calculus, and  $\DEL_L$ extends to a (necessarily unique) algebra map $\DEL_L:\Om^\bullet(M) \to G \oby  \Om^\bullet(M)$, then we say that $\Om^\bullet(M)$ is {\em covariant}. Clearly, this implies that $\Om^k \in \gmmm$, for all $k \in \bN_0$. As is easy to see, the \mpr of a covariant first-order calculus is covariant, see  \cite[\textsection 12.2.3]{KSLeabh} for details.

\subsection{Differential Calculi over $*$-Algebras}

A {\em first-order differential $*$-calculus} $(\Om^1(A),\exd)$ over a $*$-\alg $A$ is a differential calculus over $A$ \st the involution of $A$ extends to an involutive conjugate-linear map $*$ on $\Om^1(A)$ for which $(a\exd b)^*=(\exd b^*)a^*$, for all $a,b \in A$. 
If $\Om^1(G)$ is a $\ast$-calculus, then it is easy to see that the restriction of $\Om^1(G)$ to $M$ will also be a $\ast$-calculus. 

We define $\ast_\s$ to be the mapping for which the following diagram is commutative:
\bas
\xymatrix{ 
\Om^1(M)   \ar[d]_{\ast}                          & & &  \ar[lll]_{\unit \inv \circ \, (\id \,\oby \, \s \inv)}   G \coby V^1  \ar[d]^{\ast_\s} \\
\Om^1(M)       \ar[rrr]^{(\id \,\oby \, \s)\,  \circ \, \unit }                        & & & G \coby V^1.\\
}
\eas
As is routinely verified, an explicit formula for $\ast_\s$ is given by
\bal \label{astsigma}
\ast_\s \Big{(}\sum_i m^i \oby \ol{n^i}\Big{)}  =- \sum_i (m^i\1)^* \oby \ol{S (n^i)^*(m^i\2)^*}.
\eal

We call a differential calculus $\left(\Om^\bullet(A), \exd\right)$ over a $*$-\alg $A$ a {\em  $*$-differential calculus} if the involution of $A$ extends to an involutive conjugate-linear map on $\Om^\bullet$, for which $(\exd \w)^* = \exd \w^*$, for all $\w \in \Om$, and
\[
(\w_p\w_q)^*=(-1)^{pq}\w_q^*\w_p^*, \qquad \text{ for all }\w_p \in
\Om^p,~ \w_q \in \Om^q.
\]
If $\Om^1(A)$ is a first-order $*$-calculus, then its \mpr is a $*$-calculus, see  \cite[\textsection 12.2.3]{KSLeabh} for details.

\subsection{Complex Structures}

In this section we introduce a reformulation of Beggs and Smith's definition  \cite[Definition 2.6]{EBPS} of  an almost complex structure (see the remark below) which highlights alternative aspects of the structure more relevant to this and subsequent papers \cite{MMF3}. We also  recall Beggs and Smith's generalisation of  integrability to the noncommutative setting.

\begin{defn}\label{defnnccs}
An {\em almost complex structure} for a  differential $*$-calculus  $\Om^{\bullet}(A)$ over a $*$-\alg $A$ is an $\bN^2_0$-\alg grading $\bigoplus_{(a,b)\in \bN^2_0} \Om^{(a,b)}$ for $\Om^{\bullet}(A)$ such that, for all $(a,b) \in \bN^2_0$: 
\begin{enumerate}
\item \label{compt-grading}  $\Om^k(A) = \bigoplus_{a+b = k} \Om^{(a,b)}$;
\item  \label{star-cond} $\left(\Om^{(a,b)}\right)^* = \Om^{(b,a)}$.
\end{enumerate}
We call an element of $\Om^{(a,b)}$ an {\em  $(a,b)$-form}. 
\end{defn}

Let $\del$, and $\ol{\del}$ be the unique order $(1,0)$, and $(0,1)$ respectively, homogeneous operators
\bas
\del|_{\Om^{(a,b)}} : = \proj_{\Om^{(a+1,b)}} \circ \exd, & & \ol{\del}|_{\Om^{(a,b)}} : = \proj_{\Om^{(a,b+1)}} \circ \exd,
\eas
where $ \proj_{\Om^{(a+1,b)}}$, and $ \proj_{\Om^{(a,b+1)}}$, are the projections from $\Om^{a+b+1}(A)$ onto $\Om^{(a+1,b)}$, and $\Om^{(a,b+1)}$, respectively.  

\begin{lem}\label{intlem} \cite[\textsection 3.1]{EBPS}
If $\bigoplus_{(a,b)\in \bN^2_0}\Om^{(a,b)}$ is an almost-complex structure for a differential calculus $\Om^{\bullet}(A)$ over an \alg $A$, then the following conditions are equivalent:
\begin{enumerate} 
\item $\exd = \del + \ol{\del}$;
\item $\del^2=0$;
\item $\adel^2 = 0$;
\item the triple $\big(\bigoplus_{(a,b)\in \bN^2}\Om^{(a,b)}, \del,\ol{\del}\big)$ is a double complex;
\item \label{1} $\exd(\Om^{(1,0)}) \sseq \Om^{(2,0)} \oplus \Om^{(1,1)}$;
\item \label{2} $\exd(\Om^{(0,1)}) \sseq \Om^{(1,1)} \oplus \Om^{(0,2)}$.
\end{enumerate}
\end{lem}

\begin{defn}
When the conditions in Lemma \ref{intlem} hold for an almost-complex structure, then we say that it is {\em integrable}. We call an integrable almost-complex structure a {\em complex structure}, and we call the double complex  $(\bigoplus_{(a,b)\in \bN^2}\Om^{(a,b)}, \del,\ol{\del})$  the {\em Dolbeault double complex} of the complex structure.
\end{defn}

The following useful lemma shows that Beggs and Smith's definition of a complex structure is equivalent to the definition given by Khalkhali,  Landi, and  van Suijlekom in \cite[Definition 2.1]{KLVSCP1}.

\begin{lem}  \cite[\textsection 3]{EBPS}
If an almost complex structure $\bigoplus_{(a,b)\in \bN_0^2} \Om^{(a,b)}$ is integrable, then
\begin{enumerate}
\item $\del(a^*) = (\ol{\del} a)^*$, and $\ol{\del}(a^*) = (\del a)^*$,   for all $a \in A$;
\item both $\del$ and $\ol{\del}$ satisfy the graded Leibniz rule.
\end{enumerate}
\end{lem}

\begin{rem}
In \cite[Definition 2.6]{EBPS} an almost complex structure, for a differential $*$-calculus $\Om^\bullet(A)$, is defined to be a zero-order derivation $J:\Om^\bullet(A) \to \Om^\bullet(A)$ \stn, for all $a \in A$, $J(a) = 0$, and, for all $\w \in \Om^1(A)$, $J^2(\w) = -\w$ and $J(\w^*) = \left(J(\w)\right)^*$. In \cite[\textsection 2.5]{EBPS}, it is shown that, for any such $J$,  an almost complex structure in the  sense of Definition \ref{defnnccs} is uniquely determined by $J(\w) = (a-b)i\w$, for $\w \in \Om^{(a,b)}$. That the reverse construction is well-defined follows from the second part of Proposition \ref{NNC-Thm}. 
\end{rem}

\section{Quantum Projective Space}

We introduce quantum projective $N$-space $\cpn$ as a $\bC_q[U_{N-1}]$-covariant sub\alg of  the quantum special unitary group $\csun$, and describe its Heckenberger--Kolb \fodc $\Om^1_q(\ccpn)$.

\subsection{The Quantum Special Unitary Group $\csun$}

We begin by fixing notation and recalling the various definitions and constructions needed to introduce the quantum unitary group and the quantum special unitary group. (Where proofs or basic details are omitted we refer the reader to \cite[\textsection 9.2]{KSLeabh}.)

For $q \in (0,1]$ and $\nu : = q-q^{-1}$, let $\bC_q[GL_N]$ be the quotient of the free \nc \alg  ${\bC \big{<} u^i_j, \dt \inv \,|\, i,j = 1, \ldots , N \big{>}}$ by the ideal generated by the elements
\begin{align*}
u^i_ku^j_k  - qu^j_ku^i_k, &  & u^k_iu^k_j - qu^k_ju^k_i,                    & &   \; 1 \leq i<j \leq N, 1\leq k \leq N; \\
  u^i_lu^j_k - u^j_ku^i_l, &  & u^i_ku^j_l - u^j_lu^i_k-\nu u^i_lu^j_k, & &   \;  1 \leq i<j \leq N,\; 1 \leq k < l \leq N;\\
  \dt_N \dt_N \inv - 1,   & & \dt_N \inv \dt_N - 1,& &    
\end{align*}
where  $\dt_N$,  the {\em quantum determinant},  is the element
\[
\dt_{N} := \sum\nolimits_{\pi \in S_N}(-q)^{\ell(\pi)}u^1_{\pi(1)}u^2_{\pi(2)} \cdots u^N_{\pi(N)}
\]
with summation taken over all permutations $\pi$ of the set $\{1, \ldots, N\}$, and $\ell(\pi)$ is the number of inversions in $\pi$. As is well-known, $\dt_N$ is a central and grouplike element of the bialgebra.

A bi\alg structure on $\bC_q[GL_N]$ with coproduct $\DEL$, and counit $\e$, is uniquely determined by $\DEL(u^i_j) :=  \sum_{k=1}^N u^i_k \oby u^k_j$; $\DEL(\dt_N\inv) = \dt_N\inv \oby \dt_N\inv$; and $\e(u^i_j) := \d_{ij}$; $\e(\dt_N\inv) = 1$.
Moreover, we can endow  $\bC_q[GL_N]$ with a Hopf \alg structure by defining
\begin{align*}
S(\dt_{N} \inv)= \dt_{N}, ~~~~ S(u^i_j) = (-q)^{i-j}\sum\nolimits_{\pi \in S_{N-1}}(-q)^{\ell(\pi)}u^{k_1}_{\pi(l_1)}u^{k_2}_{\pi(l_2)} \cdots u^{k_{N-1}}_{\pi(l_{N-1})}\dt_N\inv,
\end{align*}
where $\{k_1, \ldots ,k_{N-1}\} = \{1, \ldots, N\}\bs \{j\}$, and $\{l_1, \ldots ,l_{N-1}\} = \{1, \ldots, N\}\bs \{i\}$ as ordered sets.
A Hopf $*$-\alg structure is determined by ${(\dt_{N}^{-1})^* = \dt_{N}}$, and $(u^i_j)^* =  S(u^j_i)$. We denote the Hopf $*$-\alg by $\bC_q[U_N]$, and call it the {\em quantum unitary group of order $N$}. We denote the Hopf  {$*$-\algn} $\bC_q[U_N]/\la \dt_{N} - 1 \ra$ by $\bC_q[SU_N]$, and call it the {\em quantum special unitary group of order $N$}.

\subsection{The Quantum Projective Spaces $\cpn$}

Following the description introduced in \cite[\textsection 3]{Mey}, we present quantum  $(N-1)$-projective space as the coinvariant sub\alg of a $\bC_q[U_{N-1}]$-coaction on $\csun$.  (This sub\alg is a $q$-deformation of the coordinate \alg of the complex manifold $SU_N/ U_{N-1}$. Recall that  $\bC P^{N-1}$ is isomorphic to $SU_N/ U_{N-1}$.)

\begin{defn}
Let $\a_N:\bC_q[SU_N] \to \bC_q[U_{N-1}]$ be the surjective Hopf \alg map  defined by setting $\a_N(u^1_1) = \dt_{N-1} \inv,$ \, $\a_N(u^1_i)=\a_N(u^i_1)=0$, for $i = 2, \cdots, N$, and  $\a_N(u^i_j) = u^{i-1}_{j-1},$ for  $i,j=2, \ldots, N$. {\em Quantum projective $(N-1)$-space} $\cpn$ is defined to be the quantum homogeneous space of the corresponding homogeneous coaction $\DEL_{SU_N,\a_N} = (\id \otimes \a_N) \circ \Delta$.
\end{defn}

As is well known,  $\cpn$ is generated as a $\bC$-\alg by the set $\{z_{ij} := u^i_1S(u^1_j)\,|\, i,j = 1,\ldots , N\}$  (see \cite[\textsection 11.6]{KSLeabh} for more details). Moreover,  $\bC_q[SU_N]$ is faithfully flat as a right module over $\cpn$ \cite{MulSch}, and so, $\cpn$ is a quantum homogeneous space. 

An important family of objects in $\gmmm$
 is the {\em quantum line bundles} $\E_p$, for $p \in \bZ$, where $\E_p := \Psi(\bC)$, with $\bC$ considered as an object in $\mhm$ 
according to the  $\bC[U_{N-1}]$-coaction  $\l \mto \l \oby \dt_{N-1}^{p}$,  for $\l \in \bC$. Clearly, we have that $\E_0 = \cpn$. (In the $q=1$ case, these modules are the modules of sections of the line bundles over $\ccpn$, see Remark \ref{rem2}.)

\subsection{The Heckenberger--Kolb Calculus $\Om^1_q(\ccpn)$}

In this subsection we recall the first-order differential calculi introduced by \hakn, and its realization as the restriction to $\cpn$ of a certain calculus on $\csun$. 

A left-covariant first-order calculus over an \alg $A$ is called {\em irreducible} if it does not possess any non-trivial quotients by a left-covariant $A$-bimodule. 
\begin{thm} \cite{HK} \label{HKClass}
There exist exactly two non-isomorphic finite-dimensional irreducible left-covariant first-order differential calculi over quantum projective  $(N-1)$-space. We call the direct sum of these two calculi the {\em Heckenberger--Kolb} calculus for $\cpn$.
\end{thm}

In general, it proves very useful to realise a calculus on a quantum homogeneous space as the restriction of a calculus over $G$.  The following proposition recalls some  details about  a calculus over $\bC_q[SU_N]$ that restricts to the \hk calculus. The technical formulae presented here will be play a crucial role in later calculations.

\begin{prop}\cite[\textsection 4, \textsection 5]{MMF1} \label{distquotient} 
For $q \neq 1$, there exists a covariant $\ast$-calculus $\Om^1_q(SU_N)$ over $\bC_q[SU_N]$ such that:
\bet

\item For $i=1, \ldots, N-1$, a basis for  $\Lambda_q(SU_N)$ is given by $e^0 :=  \ol{u^1_1-1}$, and 
\begin{flalign*}
e^+_{i} & :=  \ol{z_{i+1,1}}  = q^{-1 + \frac{2}{N}} \ol{u^{i+1}_{1}}  = - q^{-\frac{2}{N}} \ol{S(u^{i+1}_1)},\\
 e^-_{i} & := \ol{z_{1,i+1}} =  - q^{1 - 2i+ \frac{2}{N}}\ol{u^1_{i+1}}  = q^{2-\frac{2}{N}} \ol{S(u^1_{i+1})}. 
\end{flalign*}


\item For  $i,j = 2, \ldots, N$, it holds that
$
\ol{z_{ij}} =  \ol{u^i_j} = \ol{S(u^i_j)} = 0. 
$

\item \label{distquotcalc1}  For $1\leq i<j \leq N-1; k=1, \ldots, N$,  all the non-zero actions of the generators on the basis elements $e^{\pm}_i$ are given by
\begin{align} 
e^{+}_i \tl u^{j+1}_{i+1} = q^{-\frac{2}{N}}\nu e^{+}_{j-1}, & & e^-_{i} \tl u^{i+1}_{j+1} = q^{2(j-i) -\frac{2}{N}}\nu e^-_{j},
\end{align}
\begin{align}
e^{\pm}_i \tl u^{k}_{k} = q^{\d_{1k}+\d_{i+1,k}-\frac{2}{N}}e^{\pm}_i.
\end{align}

\item  For $1\leq i<j \leq N-1; k=1, \ldots, N$, all non-zero actions of the antipodes of the generators are given by
\begin{align} \label{distquotcalc2}
 e^+_i \tl S(u^{j+1}_{i+1}) = -q^{\frac{2}{N}}\nu e^+_{j}, & & e^-_i \tl S(u^{i+1}_{j+1}) = -q^{\frac{2}{N}}\nu e^-_{j}, 
\end{align}
\begin{align} \label{distquotcalc3}
e^{\pm}_i \tl S(u^k_k) = q^{\frac{2}{N} - \d_{k1}-\d_{i+1,k}}e^{\pm}_i.
\end{align}

\eet
\end{prop}

The following proposition recalls some important facts  about the restriction of $\Om^1_q(SU_N)$  to $\cpn$, principal among them that it is indeed equal to the \hk calculus.

\begin{prop}\cite[\textsection 5]{MMF1} \label{restriction} 
Denoting by $\Om^1_q(\ccpn)$ the restriction of $\Om^1_q(SU_N)$ to a  $*$-calculus over  $\cpn$, it holds that: 

\bet

\item The  right ideal  $I^{(1)} \sseq \cpn^+$, corresponding to $\Om^1_q(\ccpn)$, is generated by the elements
\begin{align} \label{icpngenset}
\{   z_{ij}, \, z_{i1}z_{j1},\,z_{1i}z_{1j}, \, z_{i1}z_{jk}, \, z_{1i}z_{jk} ~  |  ~ i,j,k = 2, \ldots, N; ~ j \neq k\},
\end{align}
which directly implies that $\Om^1_q(\ccpn)$ is an object in the subcategory $\gmm_0$.

\item A decomposition of $V^1$ in the category $\mhm$ is given by 
\bas
V^1   = & V\hol \oplus V\ahol \\
                   :  = & \spn_{\bC}\{e^+_i\,|\, i= 2, \ldots, N-1\} \oplus  \spn_{\bC}\{e^-_i\,|\, i= 2, \ldots, N-1\}.
\eas
We denote the corresponding decomposition in $\gmmm$ by 
\bas
\Om_q^1(\ccpn) := \Om_q\hol \oplus \Om_q\ahol.
\eas

\item The two calculi  $\Om_q\hol$ and $\Om_q\ahol$ are non-isomorphic, and are the calculi identified in Theorem \ref{HKClass}.

\eet
\end{prop}

\section{Monoidal Structures and Equivalences}

We use Takeuchi's categorical equivalence to transfer  the canonical monoidal structure on $\gmmm$ to a monoidal structure on $\mhm$. We then show that for the subcategory $\lmhm$ it has a particularly simple from.

\subsection{Monoidal Structures on $\gmmm$ and $\mhm$}

Let us first recall the standard monoidal structure for $\gmmm$. For $\E,\F$ two objects in $\gmmm$, we define $\E \oby_M \F$ to be the usual bimodule tensor product endowed with the standard left $G$-comodule structure
\bal \label{comodstruc1}
\DEL_{L}: \E \oby_M \F \to  G \oby \E \oby_M \F, & &  e \oby_M f \mto e \m1 f\m1 \oby e \0 \oby_M f\0.
\eal
Clearly, $\E \oby_M \F$ is well-defined as an object in $\gmmm$.
 
The equivalence between $\gmmm$ and $\mhm$ can be used to induce a monoidal structure on $\mhm$: For $V,W \in \mhm$, we define 
\bas
V \odot W := \Phi\left(\Psi(V) \oby_M \Psi(W)\right).
\eas

\subsection{The Restriction of the Monoidal Structure of $\mhm$ to the Subcategory $\lmhm$}

The explicit presentation of $\odot$ is somewhat cumbersome. However, upon restricting to the subcategory $\lmhm$ introduced in Section \ref{Takeq}, a significant simplification occurs.

This category {\em has} a natural monoidal structure $\oby$, where for $V,W$ two objects in $\lmhm$, we define $V \oby W$ to be the usual vector space tensor product, endowed with the trivial right $M$-action, and a right $H$-comodule structure given by
\bal \label{comodstruc2}
\DEL_{R}: V \oby W \to  V \oby W \oby H, & &  v \oby w \mto  v\0 \oby w\0 \oby w\1v\1.
\eal
That these two structures are compatible in the sense of (\ref{mhm-compt}) follows easily from (\ref{piofm0}). (Note that the tensor product defined here differs from the standard choice \cite[\textsection 1.3.2]{KSLeabh}.)

\begin{prop} \label{monoidal.equiv}
Denoting by $(\lmhm, \oby)$ the monoidal category for which $V \oby W$ is the standard right $H$-comodule tensor product of $V $ and $W$,
endowed with the trivial right action, a monoidal equivalence between $(\lmhm, \oby)$ and $(\lmhm, \odot)$ is given by
\bas
\mu: V \odot W \to V \oby W,  & & \ol{\big{(}\sum_i f_i \oby v_i \big{)} \oby_M \big{(}\sum_j g_j \oby w_j\big{)}} \mto \sum_{i,j} \e(f_i)\e(g_j) v_i \oby  w_j.
\eas
\end{prop}
\demo
The defining property of the subcategory $\lgmmm$ implies that an isomorphism \linebreak $\Phi(\Psi(V) \oby_M \Psi(W)) \to \Phi(\Psi(V)) \oby \Phi(\Psi(W))$ is given by
\bas
 \ol{\big{(} \sum_i f_i \oby v_i\big{)} \oby_M \big{(}\sum_j g_j \oby w_j\big{)}} \mto  \ol{\big{(}\sum_i f_i \oby v_i\big{)}} \oby \ol{\big{(}\sum_j g_j \oby w_j\big{)}}.
\eas
Composing  this isomorphism with $\frame^{\oby 2}$ gives the map $\mu$. \qed

\begin{cor}
For any covariant first-order differential calculus $\Om^1(M)$ contained in $\lgmmm$, an isomorphism is given by 
\bas
\s^{k}: \Phi\big{(}\Om^1(M)^{\oby_M k}\big{)} \to V^{\oby k},& & \ol{m_0\exd m_1 \oby \cdots \oby \exd m_k} \mto \e(m_0)\ol{m_1^+} \oby \cdots \oby \ol{m_k^+},
\eas
where $V^{\oby k} := (V^1)^{\oby k}$, and $V^1$ is the cotangent space of $\Om^1(M)$.
\end{cor}
\demo
The monoidal equivalence between $\big{(}\lmhm, \odot\big{)}$ and $\left(\lmhm, \oby\right)$ induces a unique isomorphism
$
\Phi\big{(}\Om^1(M)^{\oby_M k}\big{)} \simeq \Phi\big{(}\Om^1(M)\big{)}^{\oby k}.
$
Composing this isomorphism with $\s^{\oby k}$ gives $\s^k$ as described.
\qed

\section{Describing the  Maximal Prolongation of a Covariant First-Order Calculus}

We give explicit descriptions of the maximal prolongation of a covariant \fodcn, over a quantum homogeneous space $M$, in terms of the corresponding submodule $N \sseq \Om^1_u(M)$, and in terms of the corresponding ideal $I^{(1)} \sseq M^+$. The second presentation is then applied to the \hk calculus for $\cpn$.

Throughout this section we will assume that $\Om^1(M) \in \lgmmm$.

\subsection{Describing the  Maximal Prolongation in Terms of a Certain Submodule $I^{(2)} \sseq V^{\oby 2}$}

We show that the task of finding an explicit description of the maximal prolongation of a first-order calculus can be reduced to the problem of finding an explicit description of a certain submodule $I^{(2)} \sseq V^{\oby 2}$.

\begin{lem}\label{sigmak}
Denote  $V^k  := V^{\oby k}/I^{(k)}$, where $I^{(k)}$ is the  degree $k$ component of the ideal of  $\T(V^1)$ generated by $I^{(2)} := \s^2\left(\Phi(\exd N)\right)$. An isomorphism is given by
\bas
\s^{\wed k}:\Phi(\Om^k(M)) \to V^k, & & \ol{m_0\exd m_1 \oby \cdots \oby \exd m^k} \mto \e(m_0)\ol{m_1^+} \wed \cdots \wed \ol{m_k^+},
\eas
where we use $\wed$ to denote multiplication in $V^\bullet := \bigoplus_{k} V^k$.
\end{lem}
\demo
For $\w \in  N$, and $\exd:\Om^1_u(M) \to \Om^1(M) \oby_M \Om^1(M)$, it follows from 
\bas
\exd(m\w) = \exd m \oby \w  + m\exd \w = m\exd \w, & & \exd(\w m) = (\exd \w)m + \w \oby \exd m =  (\exd \w)m ,
\eas
that $\exd N$ is well-defined as an object in $\gmmm$. Denoting by  $\la \exd N \ra^k$ the $k\th$-component of the ideal of $\T(\Om^1(M))$ generated by  $\exd N$, exactness of $\Phi$  implies that 
\bas
\Phi\big{(}\Om^k(M)\big{)} = \P\big{(}\Om^1(M)^{\oby_M k}/\la \exd N \ra^k \big{)} \simeq \P\big{(}\Om^1(M)^{\oby_M k}\big{)}/\P\big{(}\la \exd N \ra^k\big{)}.
\eas
The map $\s^{k}$ now induces an isomorphism 
\bas
\s^k:\P\big{(}\Om^k(M)\big{)}  \simeq  V^{\oby k}/\s^{\oby k}\big{(}\P(\la \exd N \ra^k) \big{)}.
\eas
The lemma now follows from the fact that  $\s^{k}\big{(}\P(\la \exd N \ra^k)\big{)}$ is equal to the degree $k$ part of the ideal of $\T(V^1)$ generated by $I^{(2)}$.
\qed
 
The following lemma is a first step towards finding a workable description of $I^{(2)}$.
\begin{lem} \label{iideal}
It holds that
\bal \label{i2m}
I^{(2)} &   = \Big{\{}\sum_i \ol{ m_i^+} \oby \ol{n_i^+} \setl \sum_i m^i \exd n^i = 0 \Big{\}}.
\eal
\end{lem}
\demo
From the definition of the \mprn, we have
\bal \label{27}  
\Phi(\exd N)   =\Big{\{}\sum_i \ol{\exd m^i \oby_M \exd n^i} \setl  \sum_i m^i \exd n^i = 0 \Big{\}}.
\eal
Operating on (\ref{27}) by $\s^{2}$ then gives us (\ref{i2m}). 
\qed

\subsection{Describing $I^{(2)}$ in Terms of $I^{(1)}$}

While we now have an explicit description of $I^{(2)}$ in terms of $N$, it proves more useful in practice to have a description of $I^{(2)}$ in terms of $I^{(1)}$. In this section we use a certain type of first-order calculus on $G$ to produce just such a description.

\bigskip

\begin{defn}
For any \fodc $\Om^1(M)$ over $M$, a {\em framing calculus} $\Om^1(G)$ is a \fodc over $G$ such that:
\begin{enumerate}
\item   $\Om^1(M)$ is the restriction of $\Om^1(G)$ to $M$;
\item $\Om^1(M)G \sseq G\Om^1(M)$.
\end{enumerate}
\end{defn}

Some consequences of the definition of a framing calculus, under the assumption of covariance,  are presented in the following lemma.

\begin{lem}\label{framcalclem}
For $\Om^1(M)$ a covariant \fodc over $M$:

\begin{enumerate}

\item If $\Om^1(G)$ is a  covariant framing calculus, the map 
$
\iota: V^1 \to \Lambda, ~ \ol{m} \mto \ol{m}
$
is an embedding which induces a right $G$-action on $V^1$.

\item If $\Om^1(M)$ is finite dimensional, then a linear isomorphism is defined by 
\bal \label{gammadefn}
\g:V^{\oby 2} \to V^{\oby 2},  & & \ol{m^1} \oby \ol{m^2} \mto \ol{m^1 m^2\1} \oby \ol{(m^2\2)^+},
\eal
with inverse given by $\g\inv(\ol{m^1} \oby \ol{m^2}) = \ol{m^1 S(m^2\1)} \oby \ol{(m^2\2)^+}$.
\end{enumerate}
\end{lem}
\demo
The fact that $\iota$ is an embedding follows from commutativity of the  diagram
\bas
\xymatrix{ 
G \coby V^1 \ar[rrr]^{\id \oby \iota}                            & & & G \oby \Lambda \\
\Om^1(M)    \, \,\,  \ar@{^{(}->}[rrr]   \ar[u]^{\simeq}   & & & \Om^1(G) \ar[u]_{\simeq}.
}
\eas
Condition $2$ in the definition of a framing calculus implies that $\iota(V^1)$ is a $G$-submodule of $\Lambda$.

To establish the third part of the lemma, we work in the category of vector spaces. Define $\g$ to be the map for which the following diagram is commutative
\bas
\xymatrix{ 
V^{\oby 2}     \ar[d]_{(\s^{\oby 2})\inv}     \ar[rrr]^{\g}                        & & &  \Lambda^{\oby 2} \\
 \Phi\big{(}\Om^1(M)\big{)}^{\oby 2} \ar[rrr]_{E}        & & & \big{(}\Om^1(G)\big{)}^{\oby_G 2}/\big{(}G^+(\Om^1(G))^{\oby_G 2}\big{)},  \ar[u]_{\t}
}
\eas
where $E$ is the obvious linear map (well-defined since $\Om^1(M) \in \lgmmm$), and $\t$ is the isomorphism 
\bas
\t \big{(} \ol{g^0\exd g^1 \oby_G \exd g^2}\big{)} = \e(g^0)\ol{(g^1)^+g^2\1} \oby \ol{(g^2\2)^+}.
\eas
(The fact that $\t$ is  well-defined  is implied by  the canonical  isomorphism $\Om^1(G) \oby_G \Om^1(G) \simeq G \oby \Lambda^{\oby 2}$.)
To see that $\g$ coincides with the map defined in (\ref{gammadefn}) we note that
\bas
 \t \circ E \circ (\s^{\oby 2})\inv \big(\ol{m^1} \oby \ol{m^2} \big) & = \t \circ E \big(\ol{\exd m^1} \oby \ol{\exd m^2}\big) = \t  \big(\ol{\exd m^1 \oby_G \exd m^2}\big)\\
&  = \ol{m^1m^2\1} \oby \ol{(m^2\2)^+}.
\eas
Since we are assuming $V^1$ to be finite-dimensional, $\g$ is an isomorphism \iff it is surjective. This is implied by the following calculation, as is the given formula for $\g\inv$: 
\bas
\g\big(\ol{m^1S(m^2\1)} \oby \ol{(m^2\2)^+}\big) = \ol{m^1S(m^2\1)m^2\2} \oby \ol{(m^2\3)^+} =  \ol{m^1} \oby \ol{m^2}.
\eas
\qed

The following useful result is needed for the proof of the theorem below, and for the proof of Proposition (\ref{prop-integrability}). 

\begin{lem}\label{dsigma}
Let $\exd_\s$ be the unique map for which the diagram commutes
\bas
\xymatrix{ 
\Om^2(M)                        \ar[rrr]^{(\id \oby \s^{\wed 2}) \, \circ \,\unit} & & &  G \coby V^2 \\
\Om^1(M)  \ar[u]^{\exd}   & & &     \ar[lll]^{ \unit \inv \circ \, (\id \oby \s\inv)}     G\, \coby V^1 \ar[u]_{\exd_\s}.
}
\eas
For any $\sum_i f^i \oby \ol{m^i} \in G\, \coby V^1$, the map  $\exd_\s$ acts according to 
\bas
\exd_{\s}\Big(\sum_i f^i \oby \ol{m^i}\Big) = \sum_i f^i\1 \oby \ol{(f^i\2S(m^i\1))^+} \wed \ol{(m^i\2)^+} . 
\eas
\end{lem}
\demo
As a direct consequence of  Corollary \ref{unitinv} and Lemma \ref{sigmak}, it holds that
\bas
        &  (\id \oby \s^{\wed 2}) \circ \unit \circ  \exd \circ \unit \inv  \circ (\id \oby \s\inv)\Big(\sum_i f^i \oby \ol{m^i}\Big) \\
=  \,  &  (\id \oby \s^{\wed 2}) \circ  \unit  \circ \exd \Big(\sum_i f^i S(m^i\1) \exd m^i\2\Big)\\
=  \,  &   (\id \oby \s^{\wed 2}) \circ  \unit  \Big(\sum_i \exd(f^i S(m^i\1)) \wed \exd m^i\2\Big) \\
=  \, &  (\id \oby \s^{\wed 2}) \Big(\sum_i f^i\1S(m^1\2)m^i\3 \oby \ol{ \exd(f^i\2 S(m^i\1)) \wed \exd m^i\4}\Big)\\
=  \, &  (\id \oby \s^{\wed 2}) \Big(\sum_i f^i\1 \oby  \ol{\exd(f^i\2 S(m^i\1)) \wed \exd m^i\2}\Big)\\
=  \, &   \sum_i f^i\1 \oby  \ol{(f^i\2 S(m^i\1))^+} \wed \ol{(m^i\2)^+}.
\eas
\qed
An immediate consequence of the lemma is the following result.
\begin{cor}
It holds that
\bas
I^{(2)} = \Big\{ \sum_i \ol{(f^iS(z^i\1))^+} \otimes \ol{(z^i\2)^+} ~ \text{\Big{|}} ~ \sum_i f^i \oby z^i \in G\, \coby I^{(1)} \Big\}.
\eas 
\end{cor}

We now come to the central result of this section, which gives us a workable description of the ideal $I^{(2)}$.

\begin{thm}\label{theframingtheorem} 
For $\Om^1(G)$ a covariant framing calculus for $\Om^1(M)$,  we have  the equality
\bas
\iota^{\oby 2}(I^{(2)}) = \spn_{\bC}\left\{\ol{S(z\1)} \oby \ol{(z\2)^+} \,|\, z \in \text{\em Gen}(I^{(1)})\right\},
\eas
where {\em Gen}$(I^{(1)})$ is any subset of $I^{(1)}$ that generates it as a right $M$-module.
\end{thm}
\demo
It follows from the above corollary that $I^{(2)}$ is equal to
\bas
\Big\{ \sum_i \ol{(f^iS(z^i\1m\1))^+} \otimes \ol{(z^i\2m\2)^+} \,|\, m^i \in M, z^i \in \text{Gen} (I^{(1)}), \sum_i f^i \oby z^im^i \in G\, \coby I^{(1)}  \Big\}.
\eas
By applying $\iota^{\oby 2}$, and using the elementary identity $(fg)^+ = f^+g + \e(f)g^+$,  for $f,g \in G$, we see that
\bas   
          & \sum_i  \ol{(f^iS(z^i\1m^i\1))^+} \otimes \ol{(z^i\2m^i\2)^+}\\
           = &  \sum_i   \ol{(f^i)^+S(z^i\1m^i\1)} \otimes \ol{(z^i\2m^i\2)^+} +  \e(f^i)\ol{(S(z^i\1m^i\1))^+} \otimes \ol{(z^i\2m^i\2)^+}.
\eas
If we now apply $\g$, and use the elementary identity $\DEL(f^+) = f\1 \oby (f\2)^+$, for $f \in G$, we get
\bas
            &  \sum_i    \ol{(f^i)^+S(z^i\1m^i\1)z^i\2m^i\2} \otimes \ol{(z^i\3m^i\3)^+} +  \e(f^i)\ol{(S(z^i\1m^i\1))^+z^i\2m^i\2} \otimes \ol{(z^i\3m\3^i)^+}\\
           = & ~ \sum_i   \ol{(f^i)^+} \otimes \ol{(z^im^i)^+}  -  \e(f^i)\ol{(z^i\1m^i\1)^+} \otimes \ol{(z^i\2m^i\2)^+}\\
           =  & - \sum_i   \e(f^i)\ol{z^i\1m^i\1} \otimes \ol{z^i\2m^i\2} = - \sum_i   \e(f^i)\ol{z^i\1m^i\1} \otimes \ol{z^i\2}\e(m^i\2) \\
          =   &  - \sum_i    \e(f^i)\e(m^i)\ol{(z^i\1)^+} \otimes \ol{(z^i\2)^+}.
\eas
Finally, by applying $\g\inv$ we get 
\bas
          - \sum_i  \e(f^i)\e(m^i)\ol{(z^i\1)^+S(z^i\1)} \otimes \ol{(z^i\2)^+}
            & =  \sum_i \e(f^i)\e(m^i)\ol{(S(z^i\1))^+} \otimes \ol{(z^i\2)^+}\\
           &    =  \sum_i  \e(f^i)\e(m^i) \ol{S((z^i) \1)} \oby \ol{(z^i\2)^+}.
\eas
Hence, we have the inclusion
\bas
\iota^{\oby 2}(I^{(2)})  \sseq \spn_{\bC}\Big\{\ol{S(z\1)} \oby \ol{(z\2)^+} \,|\, z \in \text{\em Gen}(I^{(1)})\Big\}. 
\eas

\bigskip

For the opposite inclusion, take any $z \in I^{(1)}$, and choose a representative element in $G \coby I^{(1)}$ for the class $\counit \inv(z)$. An elementary basis argument shows that  the representative can be written in the form $1 \oby z + \sum_i g^i \oby z^i $, for some $g^i \in G^+, z^i \in I^{(1)}$. Hence, by the above calculation,  an element of $I^{(2)}$ is given by
\bas
\ol{S(z\1)} \oby \ol{(z\2)^+} + \sum_i \e(g^i) \ol{S(z^i\1)} \oby \ol{(z^i\2)^+}   = \ol{S(z\1)} \oby \ol{(z\2)^+}.
\eas
\qed

\subsection{Framing Calculi for the Heckenberger--Kolb Calculus}

We now apply the general theory developed in this section to our motivating set of examples. First, we take $\Om^1_q(SU_N)$ as a framing calculus for $\wcpn$, and use it to produce a description of $I^{(2)}$. Then we take the famous three-dimensional Woronowicz calculus $\G^1_{q}(SU_2)$ as a framing calculus for $\Omega^1_q(SU_2)$.

\subsubsection{The Calculus $\Om^1_q(SU_N)$ as a Framing Calculus for  $\wcpn$} 

In this subsection we show that  $\Om^1_q(SU_N)$ is a framing calculus for  $\Om^1_q(\ccpn)$, and use Theorem \ref{theframingtheorem} to calculate the \mpr of $\Om^1_q(\ccpn)$.

\begin{prop}\label{prop:i2:HK}
The calculus $\Om^1_q(SU_N)$ is a framing calculus for $\Om^1_q(\ccpn)$, \wrt which
the subspace $I^{(2)}$ is spanned by the elements
\bal
e^-_i \oby e^+_j + q  e^+_j \oby e^-_i, & & e^-_{i}\oby e^+_{i} + q^{2} e^+_{i}  \oby e^-_{i} + q^{-(2i+1)}\nu \sum_{a=i+1}^{N-1} q^{2a} e^+_{a} \oby e^-_a, \label{I2spanset1}\\
e^-_i \oby e^-_h + q \inv e^-_h \oby e^-_i, & & e^+_i \oby e^+_h + q  e^+_h \oby e^+_i, ~~~~~ 
e^+_i \oby e^+_i, ~~~~ e^-_i \oby e^-_i,  \label{I2spanset2}
\eal
for  $h,i,j = 1, \ldots, N-1$, $i \neq j$, and $h< i$. 
\end{prop}

\demo
To see that $\Om^1_q(SU_N)$ is a framing calculus $\wcpn$, we first recall that  $\wsun$ restricts to $\Om^1_q(\ccpn)$ on $\cpn$. The third part of Proposition \ref{distquotient} shows that $V^1$ is a right submodule of $\Lambda_{SU_N}$. Hence,  $\Om^1_q(SU_N)$  is  a framing calculus for $\wcpn$.

For sake of convenience, let us recall from Proposition \ref{restriction} that a generating set of $I^{(1)}$ is given by
\bas
\{ z_{ij},z_{i1}z_{kl}, z_{1i}z_{kl}\,|\, i,j = 2, \ldots, N, i \neq j, (k,l) \neq (1,1)\}.
\eas
For $z_{ij}$, we have 
\begin{align*}
\ol{S((z_{ij})\1)} \oby \ol{((z_{ij})\2)^+} & = \sum_{a,b=1}^N \ol{S(u^i_aS(u^b_j))} \oby \ol{(u^a_1S(u^1_b))^+} = \sum_{a,b=1}^N \ol{S^2(u^b_j)S(u^i_a)} \oby \ol{(u^a_1S(u^1_b))^+}\\
 & = \sum_{a,b=1}^N q^{2(b-j)} \ol{u^b_j S(u^i_a)} \oby \ol{(u^a_1S(u^1_b))^+}.
\end{align*}
From Proposition \ref{distquotient}, we can conclude that  the summand $\ol{u^b_j S(u^i_a)} \oby \ol{(u^a_1S(u^1_b))^+}$ is non-zero only if $a = i,b=1$, or $a=1,b=j$. Thus, we have 
\bal
\ol{S((z_{ij})\1)} \oby \ol{((z_{ij})\2)^+} & = q^{2(1-j)} \ol{u^1_j S(u^i_i)}  \oby \ol{u^i_1S(u^1_1)} +  \ol{u^j_j S(u^i_1)} \oby \ol{u^1_1S(u^1_j)} \notag \\
           & =  q^{2(1-j)}\ol{u^1_jS(u^i_i)}\oby e^+_{i-1} + \ol{S(u^i_1)u^j_j} \oby e^-_{j-1}\\
           & =  - q^{2(1- j) + 2j-3} e^-_{j-1}\oby e^+_{i-1} - e^+_{i-1}  \oby e^-_{j-1}\\
           & =  - q^{-1}e^-_{j-1}\oby e^+_{i-1} -  e^+_{i-1}  \oby e^-_{j-1},
 \label{epmrel}
\eal
where we have used Proposition \ref{distquotient}, and the standard relation $u^j_jS(u^i_1) = S(u^i_1)u^j_j$ \cite[Theorem 1]{PV}. 
This gives us the first element in (\ref{I2spanset1}). 

A similar analysis will show that the generators $z_{ii}$ give linear multiples of the element 
\bas
e^-_{i-1}\oby e^+_{i-1} + q^{2} e^+_{i-1}  \oby e^-_{i-1} + q^{1-2i}\nu \sum_{a=i}^{N-1} q^{2a} e^+_{a} \oby e^-_a,
\eas
and the generators $z_{i1}z_{i1}$, and $z_{1i}z_{1i}$ give scalar multiples of  $e^+_{i-1} \oby e^+_{i-1}$, and $e^-_{i-1} \oby e^-_{i-1}$ respectively. 

If we now assume  that $i < j$, then for $z_{i1}z_{j1}$, and $z_{j1}z_{i1}$, we have that  $\ol{S((z_{i1}z_{j1})\1)} \oby \ol{((z_{i1}z_{j1})\2)^+}$, and  $\ol{S((z_{j1}z_{i1})\1)} \oby \ol{((z_{j1}z_{i1})\2)^+}$, are both equal to linear multiples of the element 
\bas
e^+_{j-1} \oby e^+_{i-1} + q \, e^+_{i-1} \oby e^+_{j-1}.
\eas
For $z_{1i}z_{1j}$, and $z_{1j}z_{i1}$, we have that  $\ol{S((z_{1i}z_{1j})\1)} \oby \ol{((z_{1i}z_{1j})\2)^+}$, and  $\ol{S((z_{1j}z_{1i})\1)} \oby \ol{((z_{1j}z_{1i})\2)^+}$, are both equal to linear multiples of the element 
\bas
e^-_{j-1} \oby e^-_{i-1} + q \inv e^-_{i-1} \oby e^-_{j-1}.
\eas

Finally, for $k,l \neq 1$, we get that 
\bas
\ol{S((z_{i1}z_{kl})\1)} \oby \ol{((z_{i1}z_{kl})\2)^+} = \ol{S((z_{1i}z_{kl})\1)} \oby \ol{((z_{1i}z_{kl})\2)^+} = 0.
\eas

\qed

\subsubsection{The Woronowicz $3D$-Calculus on $\wsu2$ as a Framing Calculus}

In this subsection we specialise to the case of $\cp1$, and use the famous Woronowicz $3D$-calculus $\G^1_{q}(SU_2)$ over $\bC_q[SU_2]$ as a framing calculus. This serves to highlight the fact that there can exist more than one framing calculus for a calculus over a quantum homogeneous space. 

Following standard convention,  for the special case of $\csu2$, we use a simplified notation for the generators: $a := u^1_1, b := u^1_2, c := u^2_1$, and $d: = u^2_2$. Moreover, since $V\hol$ and $V\ahol$ are both $1$-dimensional, we will denote $e^+_1$ by $e^+$, and $e^-_1$ by $e^-_1$.
\begin{defn}
The {\em Woronowicz $3D$-calculus} $\G^1_q(SU_2)$ is the left-covariant \fodc over $\csu2$ corresponding to the ideal of $\csu2^+$ 
generated by the elements 
\bal \label{Woro3drightrels}
a+q^{-2}d-(1+q^{-2})1,\, bc, \, b^2, \, c^2, \, (a-1)b,\, (a-1)c.
\eal
\end{defn}

Using $\G^1_q(SU_2)$ as a framing calculus, we calculate $I^{(2)}$, and see that it agrees with Proposition \ref{prop:i2:HK} for the case of $N=2$.

\begin{lem}
The calculus $\G^1_q(SU_2)$ is a framing calculus for $\Om^1_q(\ccp1)$, and 
\bal \label{I2cp1}
I^{(2)} = \spn_{\bC}\{e^+ \oby e^+, \, e^- \oby e^-, \, e^- \oby e^+ + q^{2} e^+ \oby e^-\}.
\eal
\end{lem}
\demo
For the case of $N=2$,  the description of the  ideal corresponding to the calculus given in the second part of Proposition \ref{distquotient} reduces to $I^{(1)}= \la (ab)^2,bc,(cd)^2\ra$. Hence, we have a well-defined map $V^1 \to \Lambda_{SU_2}$. We know from  \cite[Example 2]{Wor} that the elements $\ol{b}$ and $\ol{c}$ are linearly independent in $\Lambda_q(SU_2)$. From (\ref{Woro3drightrels}) we have that 
$e^+ =  \ol{cd} = q \ol{c}$, and $e^- = -q \inv \ol{ab} = - q\inv \ol{b}$, and so, the map is an inclusion. Moreover, since it is also clear from (\ref{Woro3drightrels}) that $V^1$ is a right $\bC_q[SU_2]$-submodule of $\Lambda_{SU_2}$,  we have that $\G^1_q[SU_2]$ is a framing calculus for $\Om^1_q(\ccp1)$.

We now come to the calculation of $I^{(2)}$: For $bc$, we have 
\bas
\ol{S((bc)\1)}\oby \ol{((bc)\2)^+} & = \ol{S(ac)} \oby \ol{ba} + \ol{S(ad)} \oby \ol{bc} + \ol{S(bc)} \oby \ol{(da)^+} + \ol{S(bd)} \oby \ol{dc}\\
    & = -q\, \ol{cd} \oby \ol{ba} - q\inv \ol{ab} \oby \ol{dc} = qe^+ \oby e^- + q\inv e^- \oby e^+\\
& = q\inv(e^- \oby e^+ + q^{2}e^+ \oby e^-).
\eas
Analogous calculations will show that $\ol{S(((ab)^2)\1)} \oby \ol{(((ab)^2)\2)^+}$ is equal to a scalar multiple of $e^- \oby e^-$, and $\ol{S(((cd)^2)\1)} \oby \ol{(((cd)^2)\2)^+}$ is equal to a scalar multiple of $e^- \oby e^-$.
\qed

\section{Covariant Almost Complex Structures}

We begin this section by producing sufficient and necessary conditions for a bimodule decomposition of a \fodc to extend to a  complex structure on its maximal prolongation. 
The notion of a factorisable complex structure is then introduced, and a convenient formulation of the concept at the level of $1$-forms is given.
Finally, we apply our results to the  Heckenberger--Kolb calculi for the quantum projective spaces.

\subsection{Extending $1$-Form Decompositions to Almost Complex Structures}

For any smooth manifold $M$, every decomposition of the cotangent bundle into a direct sum of sub-bundles of equal dimension extends to an almost complex structure on $M$. As the following proposition shows, things are more complicated in the \nc setting. 

\begin{prop} \label{NNC-Thm}
If $\Om\hol$ and  $\Om\ahol$ are sub-bimodules of $\Om^1(A)$ \st $\Om^1(A) = \Om\hol \oplus \, \Om\ahol$:
\begin{enumerate}
\item the decomposition has at most one extension, satisfying part 1 of Definition \ref{defnnccs}, to an $\bN^2_0$-grading of the \mpr of $\Om^1(A)$;

\item such an extension exists \iff $\exd N$ is homogeneous \wrt the decomposition 
\bal \label{2:Decomp:rest}
(\Om^1(A))^{\oby_A 2} \simeq \T_{\Om}^{(2,0)} \oplus \T_{\Om}^{(1,1)} \oplus \T_{\Om}^{(0,2)},
\eal
where $\T_{\Om}^{(\bullet,\bullet)}$ denotes the unique $\bN^2_0$-grading of $\T(\Om^1(A))$ extending the decomposition of $\Om^1(A)$.

\item Moreover, condition 2 holds  \iff $(\Om^{(1,0)})^* = \Om^{(0,1)}$, or equivalently \iff $(\Om^{(0,1)})^* = \Om^{(1,0)}$.
\end{enumerate}
\end{prop}

\demo
Since $\la \exd N \ra$ is generated as an ideal by $\exd N$, it is clear that homogeneity of $\exd N$, \wrt the decomposition in  (\ref{2:Decomp:rest}), will imply homogeneity of $\la \exd N \ra$ \wrt the grading $\T_{\Om}^{(\bullet,\bullet)}$. Hence, $\T_{\Om}^{(\bullet,\bullet)}$ will descend to a grading on the maximal prolongation.  Conversely, if  $\exd N$ is not homogeneous \wrt to (\ref{2:Decomp:rest}), then $\T_{\Om}^{(\bullet,\bullet)}$ obviously  cannot descend to a grading on the \mprn.

Next we  show that this grading is the only possible $\bN^2_0$-grading on the \mpr extending the decomposition of $\Om^1(A)$: For another distinct grading $\Gamma^{(\bullet,\bullet)}$ to exist, there would have to be an element $\w  \in \T_{\Om}^{(a,b)}$, for some $(a,b) \in \bN^2_0$, such that the image of $\w$ in $\Om^{\bullet}(A)$ was not contained in $\Gamma^{(a,b)}$. Now every element of $\T_{\Om}^{(a,b)}$ is of the form
$
\w := \sum_{i=1} \w^i_1 \oby \cdots \oby  \w^i_{a+b},
$
where each summand $ \w^i_1 \oby \cdots \oby  \w^i_{a+b}$ has exactly $a$ of its factors contained in $\Om\hol$, and  $b$ of its factors  contained in $\Om\ahol$. The general properties of a graded \alg imply that the image of such an element in $\Om^{\bullet}(A)$ must be contained in $\G^{(a,b)}$. Hence,  we must conclude that there exists no other grading on the \mpr extending the decomposition $\Om^1(A)$. This gives us the first and second parts of the proposition.

\bigskip

For the third and final part of the theorem,  note that  since the $*$-map is involutive,  assuming  $(\Om^{(1,0)})^* = \Om^{(0,1)}$ is equivalent to assuming $(\Om\ahol)^* = \Om\hol$.  Every element of $\Om^{(a,b)}$ is of the form $\w := \sum_{i=1} \w^i_1 \wed \cdots \wed  \w^i_{a+b}$, where each summand $ \w^i_1 \wed \cdots \wed  \w^i_{a+b}$ has exactly $a$ of its factors contained in $\Om\hol$, and  $b$ of its factors  contained in $\Om\ahol$.  The properties of a graded $*$-\alg imply that 
\bal
\w^* = \sum_{i=1} (\w^i_1 \wed \cdots \wed  \w^i_{a+b})^* = \sum_{i=1} (-1)^{\frac{(a+b)(a+b-1)}{2})} (\w^i_{a+b})^* \wed \cdots \wed  (\w^i_{1})^*.
\eal
Our two equivalent assumptions, and the properties of a graded \algn, now imply that $\w^*$ must be contained in $\Om^{(b,a)}$, giving us that  $(\Om^{(a,b)})^* \sseq \Om^{(b,a)}$. The opposite inclusion is established analogously, giving us the desired equality.
\qed

\begin{rem}
The first part of Proposition \ref{NNC-Thm} is a special case of a  general result about graded modules \cite[II \textsection 5.5]{BourAlg}. Parts 2 and 3 can also be proved at this level of generality. We do not do so here so as to avoid unnecessary abstraction.
\end{rem}

\subsection{Covariant Almost Complex Structures}

In this subsection, we introduce the notion of a covariant almost complex structure, and find a set of simple conditions for such a structure to exist.

\begin{defn}
An almost complex structure $\Om^{(\bullet,\bullet)}$ for a covariant differential $*$-calculus $\Om^\bullet(M)$ is {\em left-covariant} if $\Om^{(a,b)}$ is a sub-object of $\Om^k$ in $\gmmm$, for all  $(a,b) \in \bN^2_0$. 
\end{defn}

The following theorem is an easy consequence of Proposition \ref{NNC-Thm} and Proposition \ref{monoidal.equiv}, so we omit the proof.

\begin{thm} \label{cor:complexstructure}
For a covariant differential $*$-calculus $\Om^\bullet(M)$:

\begin{enumerate}

\item An almost complex structure is covariant \iff  $\Om^{(1,0)}$ and $\Om^{(0,1)}$ are objects in $\gmmm$.


\item
If $\Om^1(M) \in \lgmmm$, then such a decomposition extends to an $\bN^2_0$-grading of the \mpr of $\Om^1(M)$ \iff $I^{(2)}$ is homogeneous \wrt the decomposition 
\bal \label{vmdecomposition}
V^{\oby 2} = T_V^{(2,0)} \oplus \T_V^{(1,1)} \oplus \T_V^{(0,2)},
\eal
where  $\T_V^{(\bullet,\bullet)}$ is the obvious grading of $\T(V^1)$ induced the isomorphisms $\s^k$.

\end{enumerate}
\end{thm}


We finish this subsection by deriving a simple pair of sufficient conditions for the second axiom of an almost complex structure to hold.

\begin{prop}\label{cov:star:cond}
Let $\Om^1(M) = \Om\hol \oplus \Om\ahol$ be a decomposition in $\gmmm$, and let $\Om^1(G)$ be a covariant $\ast$-calculus on $G$ that frames $\Om^1(M)$. If
\bal \label{respect}
\Om\hol G \sseq G\Om\hol, & & \text{ and } & & \Om\ahol G \sseq G\Om\ahol,
\eal
then we have $(\Om^{(1,0)})^* = \Om^{(0,1)}$ \iff
\begin{align} \label{SofV}
\big\{\ol{S(m)^*} \,|\, \ol{m} \in V\hol \big\} = V\ahol, && \text{ or equivalently }&&\big\{\ol{S(m)^*} \,|\, \ol{m} \in V\ahol\big\} = V\hol.
\end{align}
\end{prop}

\demo
If the first inclusion in (\ref{respect}) holds, then $V\hol$ is a right $G$-submodule of $\Lambda$. Taken together with (\ref{astsigma}), this then implies that $(\Om^{(1,0)})^* \sseq \Om^{(0,1)}$.  That the second equality in  (\ref{SofV}) is equivalent to the first follows from the identity $S(S(g)^*)^* = g$, for all $g \in G$. This similarly implies that $(\Om\ahol)^* \sseq \Om\hol$, giving the required equality $(\Om^{(1,0)})^* = \Om^{(0,1)}$.

Conversely, working in the category $\ggmg$, it is easy to see that $(\id \oby \s) \circ \unit: \Om^1(G) \to G \oby \Lambda$ restricts to  isomorphisms $G\Om\hol \simeq G \oby V\hol$, and $G\Om\ahol \simeq G \oby V\ahol$, where as usual $\Lambda$ denotes the space of left-invariant $1$-forms of $\Om^1(G)$. Hence, we have the commutative diagram 
\bas
\xymatrix{ 
G \oby V\hol   \ar[rrrr]^{\simeq}                     \ar[d]_{\ast_\s} & & &  & G \Om\hol \ar[d]^{\ast} \\
~~~~~~~~ ~~~~ ~~~~~G \oby V^1 \supseteq G \oby V\ahol     & & &   &    G \Om\ahol  \supseteq \Om\ahol G    \ar[llll]^{\simeq}.  ~~~~~~~~~~~~~ 
}
\eas
Hence, $\ast_\s(1 \oby \ol{m}) = 1 \oby \ol{S(m^*)} \in G \oby V\ahol$ giving $\{\ol{S(m)^*} \,|\, \ol{m} \in V\hol\} \sseq V\ahol$. 
The inclusion $\{\ol{S(m)^*} \,|\, \ol{m} \in V\ahol\} \sseq V\hol$ is established analogously. This  implies that
\bas
V\hol = \big\{\ol{S(S(m)^*)^*} \,|\, \ol{m} \in V\hol\big\} \sseq \big\{\ol{S(m)^*} \,|\, \ol{m} \in V\ahol\big\} \sseq V\hol,
\eas  
giving the first equality in  (\ref{SofV}). The second equality is established similarly.
\qed

\subsection{Factorisable Almost Complex Structures}

In this section we introduce the property of factorisability for an almost complex structure. The Dolbeault double complex of every complex manifold automatically satisfies this property \cite[\textsection 1.2]{HUY}, as do the \hk calculi for the all irreducible flag manifolds \cite[Proposition 3.11]{HKdR}.

\begin{defn}
An {\em almost complex structure} for a  differential $*$-calculus  $\Om^{\bullet}(A)$ over a $*$-\alg $A$, is called {\em factorisable} if we have bimodule isomorphisms
\bal  \label{wedge-cond} 
\wed:\Om^{(a,0)} \oby_A \Om^{(0,b)} \simeq \Om^{(a,b)},  & & \text{ and } & & \wed: \Om^{(0,b)} \oby_A \Om^{(a,0)} \simeq \Om^{(a,b)}.
\eal
\end{defn}

The following proposition  establishes a simple set of necessary and sufficient criteria for an \acs to be factorisable. 
 
\begin{prop}
An almost complex structure is factorisable \iff  we have bimodule isomorphisms 
\bal \label{wedge-cond-reform}
\wed:\Om^{(1,0)} \oby_A \Om^{(0,1)} \to \Om^{(1,1)}, & & \wed: \Om^{(0,1)} \oby_A \Om^{(1,0)} \to \Om^{(1,1)}.
\eal
\end{prop}
\demo
Surjectivity of the first map in (\ref{wedge-cond-reform}) means that for any $\w^+ \in \Om \hol$, and $\w^- \in \Om \ahol$, there exist forms $\w_i^+ \in \Om\hol$, and $\w^-_i \in \Om\ahol$ \st $\w^- \wed \w^+ = \sum_i \w^+_i \wed \w^-_i$. This easily implies surjectivity of the first map in (\ref{wedge-cond}). The proof in the other direction is trivial. That surjectivity of the second map in (\ref{wedge-cond}) is equivalent to surjectivity of the second map in (\ref{wedge-cond-reform}) is established  analogously.

The first map in (\ref{wedge-cond}) is injective if, for all $(a,b) \in \bN^2_0$, 
\bal \label{subspace}
\la \exd N \ra \cap (\T_{\Om}^{(a,0)} \oby_A \T_{\Om}^{(0,b)}) =  \la \exd N\ra^{(a,0)} \oby_A \T_{\Om}^{(0,b)} + \T_{\Om}^{\oby (a,0)} \oby_A \la \exd N\ra^{(0,b)}, 
\eal
where $ \la \exd N\ra^{(a,0)}$, and $ \la \exd N\ra^{(0,b)}$, are respectively the $(a,0)$, and $(0,b)$, homogeneous components of $\la \exd N \ra$, \wrt the grading $\T_{\Om}^{(\bullet,\bullet)}$. Now when the first mapping in (\ref{wedge-cond-reform}) is injective 
\bas
\exd N \cap \big(\Om \hol \oby_A \Om \ahol\big) =\{0\}.
\eas 
Thus, for a general element $\sum_i \nu_i \oby \w_i \oby \nu'_i \in \la \exd (N) \ra$,  where $\nu_i,\nu_i' \in \T(\Om^1(A))$, and  $w_i \in \exd N$ is a homogeneous element of $\exd N$, when  $\sum_i \nu_i \oby \w_i \oby \nu'_i \in \T_{\Om}^{(a,0)} \oby_A \T_{\Om}^{(0,b)}$ we must have that  $\w^i \in \T_{\Om}^{(2,0)}$, or $\w^i \in \T_{\Om}^{(0,2)}$. Hence (\ref{subspace}) holds, and  the first map in (\ref{wedge-cond-reform}) is injective. The proof in the other direction is trivial. That injectivity of the second map in (\ref{wedge-cond}) is equivalent to injectivity of the second map in (\ref{wedge-cond-reform}) is established  analogously. 
\qed

Finally, we specialise to the case of a covariant almost complex structure, such that $\Om^1(M) \in \lgmmm$, and find a very useful reformulation of (\ref{wedge-cond-reform}). We omit the proof which follows directly from the above proposition and Proposition \ref{monoidal.equiv}.

\begin{cor} \label{67}
If $\Om^1(M) \in \lgmmm$, then  the almost complex structure is factorisable \iff
 isomorphisms in $\lmhm$ are given by
\bal \label{wedsigma}
\wed:V \hol \oby V \ahol \to V^{(1,1)}, & &  \wed:V \ahol \oby V \hol \to V^{(1,1)}.
\eal
\end{cor}

\subsection{A Factorisable Almost Complex Structure for the Maximal Prolongation of the Heckenberger--Kolb Calculus}

We now apply the general theory developed in this section to the \hk calculus. (Note that the isomorphisms in (\ref{bundleiso}) are direct generalisations of well known classical results \cite[\textsection 2.4]{HUY}. In particular, the third isomorphism generalises orientability of $\ccpn$.)

\begin{prop} 
The decomposition  $\Om^1_q(\ccpn) = \Om\hol \oplus \Om \ahol$ (as presented in Proposition \ref{distquotient}) 
extends to a covariant factorisable almost complex structure on the \mpr of $\Om^1_q(\ccpn)$. 
\end{prop}
\demo
That $I^{(2)}$ is homogeneous \wrt the decomposition (\ref{vmdecomposition}) in Theorem \ref{cor:complexstructure} follows directly from Proposition \ref{prop:i2:HK}, as does the fact that the maps (\ref{wedsigma}) in Corollary \ref{67} are isomorphisms.  The inclusions (\ref{respect}) in Proposition \ref{cov:star:cond} follow directly from the third part of Proposition \ref{distquotient}. Moreover,  (\ref{SofV}) follows from the fact that, for $i = 2, \ldots, N$, we have
\bas
\ol{S(z_{i1})^*} &  = \ol{S(u^i_1S(u^1_1))^*} = \ol{(u^1_1S(u^{i}_1))^*} = \ol{S(u^i_1)^*(u^1_1)^*}=  \ol{S \inv((u^i_1)^*)S(u^1_1)} \\
                           &  = \ol{S \inv \circ S (u^1_i)S(u^1_1)}   =\ol{u^1_iS(u^1_1)} = q^{-4+2i} e^-_{i-1} \in V\ahol,
\eas
where we have used the standard Hopf $*$-\alg identity $\ast \circ S = S\inv \circ \ast$.
Thus, the decomposition  extends to an almost complex structure.
\qed

\begin{cor}
The vector space dimension of $V^{(a,b)}$ is $\binom{N-1}{a}\binom{N-1}{b}$, and a basis is given  by
\begin{align*}
\{ e^+_{i_1} \wed \cdots \wed e^+_{i_a} \wed e^-_{j_1} \wed \cdots \wed e^-_{j_b}\, | \, i_1 < \cdots < i_a; \, j_1 < \cdots < j_b \}.
\end{align*}
\end{cor}
\demo
That the proposed basis spans $V^{\bullet}$ is obvious from the set of generators of  $I^{(2)}$ given in Proposition \ref{prop:i2:HK}. Hence, we only need  to establish  linear independence. We begin with the case of $V^{(N-1,0)}$, where this amounts to showing that we have a non-zero vector space. To this end, we define a function $f:(V\hol)^{\oby N -1} \to \bC$ by specifying its values on the basis $\{e^+_{i_1} \oby \cdots \oby e^+_{i_{N-1}} \,|\,   i_1, \ldots, i_{N-1} = 1, \ldots, N-1\}$.
For basis elements of the form $e^+_{\pi(1)} \oby \cdots \oby e^+_{\pi(N-1)}$, where $\pi$ is a permutation of the set $\{1, \ldots, N-1\}$, we define 
\bas 
f( e^+_{\pi(1)} \oby \cdots \oby e^+_{\pi(N-1)}) : = (-q)^{-\ell(\pi)},
\eas
where $\ell$ is the number of inversions in $\pi$. On all other basis elements we set $f$ to zero. It follows from the description of the generators of $I^{(2)}$ given in Theorem \ref{prop:i2:HK}, and the definition of $I^{(N-1)}$, that $f$ descends to a non-zero map on the quotient $V^{(N-1,0)} = (V\hol)^{\oby N-1}/I^{(N-1)}$. Hence,  $V^{(N-1,0)} \neq 0$. 

We now move on to the case of $V^{(a,0)}$, for $a = 2, \ldots, N-2$. Suppose we have a linear combination of the proposed basis vectors, with each summand of degree $(a,0)$, and for which
\bas
\sum_{i_1 < \cdots < i_a} \l_{i_1, \ldots,  i_a} e^+_{i_1} \wed \cdots \wed e^+_{i_a} = 0.
\eas
Now for any $1 < j_1 \leq \cdots \leq j_a \leq N-1$, if we denote by  $\{j'_1, \ldots, j'_{N-a-1}\}$ the set complement  of $\{j_1, \ldots, j_a\}$ in  $\{1, \ldots, N-1\}$, then
\bas
 \bigg{(}\sum_{i_1 < \cdots < i_a} \l_{i_1, \ldots,  i_a} e^+_{i_1} \wed \cdots \wed e^+_{i_a}\bigg)  \wed (e^+_{j'_1} \wed \cdots  \wed e^+_{j'_{N-a-1}}{)}  \\
 = \, \l_{j_1, \ldots,  j_a} e^+_{j_1} \wed \cdots \wed e^+_{j'_a}  \wed e^+_{j'_1} \wed \cdots \wed e^+_{j'_{N-a-1}}\\
=  \,\l_{j_1, \ldots,  j_a} (-q)^m e^+_1 \wed \cdots \wed e^+_{N-1}, ~~~~~~~~~~~~~~~
\eas
for some $m \in \bZ$.  Since $e^+_1 \wed \cdots \wed e^+_{N-1} \neq 0$, we must have $\l_{j_1, \ldots,  j_a} = 0$. Thus, the proposed basis elements which are contained in $V^{(\bullet, 0)}$ are linearly independent. 

A similar argument will establish linear independence of the basis elements contained in $V^{(0,\bullet)}$. Linear independence of all the proposed basis elements now follows from the fact that the calculus is factorisable.
\qed

The above corollary tells us that the bundles $\Om^{(N-1,0)}$, $\Om^{(0,N-1)}$, and $\Om^{(N-1,N-1)}$ are non-trivial line bundles. The lemma below identifies these line bundles in terms of the integer classification of equivariant line bundles.

\begin{lem}
The right $\bC_q[U_{N-1}]$-coactions  on  $V^{(1,0)}$ and $V^{(0,1)}$ are given by
\bas
\DEL_{R}(e^+_{i}) =  \sum_{k=1}^{N-1} e^+_{k} \oby  S(u^{i}_{k})\dt_{N-1} \inv,& & \DEL_R(e^-_{i}) =  \sum_{k=1}^{N-1} e^-_{k} \oby  S^2(u^{k}_{i})\dt_{N-1}. 
\eas
Moreover, 
\begin{align} \label{bundleiso}
\Om^{(N-1,0)} \simeq \E_{-N},  & & \Om^{(0,N-1)} \simeq \E_{N}, & &  \Om^{(N-1,N-1)} \simeq  \cpn.
\end{align}
 \end{lem}
\demo
We begin by calculating the right $\bC_q[U_{N-1}]$-coaction on the basis elements of $V^{(1,0)}$. For $e^+_i$, we have
\begin{align*}
\DEL_R(e^+_{i}) = \DEL_R(\ol{z_{i+1,1}}) & = \sum_{k,l=1}^N \ol{u^k_1S(u^1_l)} \oby  S(\pi(u^{i+1}_k)S(\pi(u^l_1)))\\
                                                     & = \sum_{k=2}^{N} \ol{u^k_1S(u^1_1)} \oby  S^2(\pi(u^1_1))S(u^{i}_{k-1})\\
                                                     & = \sum_{k=2}^{N} \ol{z_{k1}} \oby \dt_{N-1} \inv  S(u^{i}_{k-1})\\
                                                     & = \sum_{k=1}^{N-1} e^+_{k} \oby  S(u^{i}_{k})\dt_{N-1} \inv.
\end{align*}
Hence, the coaction on $V^{(N-1,0)} \simeq \bC e^+_1 \wed \cdots \wed e^+_{N-1}$, which we denote by $\DEL^{N-1}_R$, acts according to
\begin{align*}
\DEL^{N-1}_R (e^+_1 \wed \cdots \wed e^+_{N-1}) & = \sum_{l=1}^{N-1} \sum_{k_l =1}^{N-1}  e^+_{k_1} \wed \cdots \wed e^+_{k_{N-1}} \oby                                                                            S(u_{k_{N-1}}^{N-1}) \cdots  S(u^1_{k_1})\dt^{1-N}_{N-1}\\
                                                & = \sum_{l=1}^{N-1} \sum_{k_{l}=1}^{N-1}  e^+_{k_1} \wed \cdots \wed e^+_{k_{N-1}} \oby
                                                     S(u^{1}_{k_1} \cdots  u^{N-1}_{k_{N-1}})\dt^{1-N}_{N-1}.
\end{align*}
Since any summand with a repeated basis element in the first tensor factor will be zero,
\begin{align*}
 \DEL^{N-1}_R(e^+_1 \wed \cdots e^+_{N-1}) & = \sum_{\pi \in S_{N-1}}  e^+_{\pi(1)} \wed \cdots \wed e^+_{\pi(N-1)} \oby
                                               S(u^{1}_{\pi(1)} \cdots  u^{N-1}_{\pi(N-1)})\dt^{1-N}_{N-1}.
\end{align*}
Moreover,  $e^+_{\pi(1)} \wed \cdots \wed e^+_{\pi(N-1)}  = (-q)^{\ell(\pi)} e^+_{1} \wed \cdots \wed e^+_{N-1}$, for any $\pi \in S_{N-1}$. Since
\begin{align*}
\sum_{\pi \in S_{N-1}}  (-q)^{\ell(\pi)}u^{1}_{\pi(1)}  \cdots  u^{N-1}_{\pi(N-1)}= \dt_{N-1},
\end{align*}
we must have  
\begin{align*}
\DEL^{N-1} (e^+_1 \wed \cdots \wed e^+_{N-1}) & = e^+_{1} \wed \cdots \wed e^+_{N-1} \oby S(\dt_{N-1})\dt_{N-1}^{1-N}\\
                                           & = e^+_{1} \wed \cdots \wed e^+_{N-1} \oby \dt^{-N}_{N-1},
\end{align*}
which in turn implies that  $\Om^{(N-1,0)} \simeq \E_{-N}$.

An analogous argument will establish that $\Om^{(0,N-1)}$ is isomorphic to $\E_{N}$. It follows as a direct consequence of these two results that $\Om^{(N-1,N-1)}$ is isomorphic to $\cpn$.
\qed

\section{Complex Structures}

In this section we give a simple set of sufficient criteria for a covariant almost complex structure to be integrable, find an interesting connection between integrability and the maximal prolongations of $\Om\ahol$ and $\Om\ahol$, and show that the Heckenberger--Kolb calculus for $\cpn$ satisfies these criteria.



\subsection{Integrability for a Covariant Almost Complex Structure}

We use the assumption of covariance to find a simple set of sufficient criteria for an almost-complex structure to be integrable. Throughout this subsection $\Om^{(\bullet,\bullet)}$ 
denotes a covariant almost-complex structure, \st  each $\Om^{(a,b)}$ is an object in $\lgmmm$.

\begin{prop} \label{prop-integrability}
The almost-complex structure $\Om^{(\bullet,\bullet)}$ is integrable if, for any linear projection $P: \Lambda^{\oby 2} \to V^{\oby 2}$, it holds that
\bal \label{int:1}
P\big(\ol{(S(m\1))^+} \oby \ol{m\2^+}\big) \in \iota^{\oby 2} \big(\T_V^{(2,0)}\big) \oplus \iota^{\oby 2}\big(\T_V^{(1,1)}\big), & & \text{for all\,\, } \ol{m} \in V\hol,
\eal
where $\iota: V^{1} \to \Lambda$ is the embedding introduced in Lemma \ref{framcalclem}. Integrability also follows from the corresponding condition for  $V\ahol$.
\end{prop}
\demo
Denote by $K:V^2 \to \Lambda^{\oby 2}/I^{(2)}$ the obvious linear embedding, and $\ol{v} \wed \ol{w} := \proj(\ol{v} \oby \ol{w})$, where $\proj: \Lambda^{\oby 2} \to \Lambda^{\oby 2}/I^{(2)}$ is the canonical projection. By  Lemma \ref{sigmak}, for any $\sum_i f^i \oby \ol{m^i} \in G \, \coby V\hol$, we have that 
\bas
   &(\id \oby K) \circ \exd_\s \Big{(}\sum_i f^i \oby \ol{m^i}\Big{)}\\
= &   \, (\id \oby K) \Big{(}\sum_i f^i\1 \oby \ol{(f^i\2S(m^i\1))^+} \wed \ol{(m^i\2)^+}\Big{)} \\
= &   \sum_i f^i\1 \oby \e(f^i\2)\ol{\r(S(m^i\1))^+} \wed \ol{(m^i\2)^+} +  \sum_i f^i\1 \oby \ol{(f^i\2)^+S(m^i\1)} \wed \ol{(m^i\2)^+} \\
= &   \sum_i f^i \oby \ol{(S(m^i\1))^+} \wed \ol{(m^i\2)^+}  + \sum_i f^i\1 \oby \ol{(f^i\2)^+\e(S(m^i\1))} \wed \ol{(m^i\2)^+}\\
= &   \sum_i f^i \oby \ol{(S(m^i\1))^+} \wed \ol{(m^i\2)^+}  + \sum_i f^i\1 \oby \ol{(f^i\2)^+} \wed \ol{(m^i)^+},
\eas
where in the third line we have used the standard identity $(fg)^+ = \e(f)g^+ + f^+g$, for $f,g \in G$, and in the fourth line we have used the assumption that $\Om^{(a,b)} \in \lgmmm$, for all $(a,b) \in \bN^2_0$. 
Thus, if (\ref{int:1}) holds then $\exd_\s \big(\sum_i f^i \oby \ol{m^i}\big) \in G \, \coby \big(V^{(2,0)} \oplus V^{(1,1)} \big)$, which in turn implies that $\exd(\Om^{(1,0)}) \sseq \Om^{(2,0)} \oplus \Om^{(1,1)}$.
 
That integrability follows from the corresponding condition for  $V\ahol$ is established analogously.
\qed

\subsection{Integrability and the Maximal Prolongations of $\Om\hol$ and $\Om\ahol$}

In this subsection we give an alternative characterisation of integrability for a (not necessarily covariant) almost complex structure $\Om^{(\bullet,\bullet)}$:The pairs $(\Om\hol,\del)$ and $(\Om\ahol,\adel)$ are each first-order differential calculi with their own respective maximal prolongations. Let us denote the $k$-forms of the \mpr of $\Om\hol$, and $\Om\ahol$, by $(\Om\hol)^k$, and  $(\Om\ahol)^k$ respectively.  

\begin{lem} \label{lem:altconst}
For an almost complex structure $\Om^{(\bullet,\bullet)}$,  the equalities 
\bal \label{intequivcond}
\big(\Om\hol\big)^k  = \Om^{(k,0)}, & & \text{ and } & & \big(\Om\ahol\big)^k = \Om^{(0,k)},
\eal
are equivalent to each other, and to integrability.
\end{lem}
\demo
Let  $\{\w_i^- \}_i$, be a subset of $\Om^1_u(M)$, \st 
$ 
 \spn_{\bC}\{\w_i^-\} = \Om\ahol,
$
(where  we  use the same symbol for $\w^{-}_i$ as for its coset in $\Om^1(A)$). If $N$ is the sub-bimodule of $\Om^1_u(A)$ corresponding to $\Om^1(A)$, then it is clear that the sub-bimodule of $\Om^1_u(A)$ corresponding to   $\left(\Om\hol,\del \right)$ is given by 
\bas
N^+ := N \, + \, \spn_{\bC}\{\w^-_i\}_i. 
\eas 
Since $(\Om\hol)^{\otimes_A k} = \T_{\Om^1}^{(k,0)}$, it is clear that the lemma would follow from the equality $\del N^+ = \la\exd N\ra^{(2,0)}$ (where by abuse of notation we mean $\exd$ and $\del$ in the sense of (\ref{abused})).
But $\del B = \la\exd B\ra^{\oby (2,0)}$, for any bimodule $B \sseq \Om^1_u(A)$. Hence,  having $\del N^+ = (\exd N)^{\oby (2,0)}$  is equivalent to having $\del \w^-_i = 0$, for all $i$, which is in turn equivalent to the almost complex structure being integrable.

That the second  equality in (\ref{intequivcond}) is equivalent to integrability is established analogously.
\qed

\subsection{Integrability of the \hk Calculus}

We now apply the general results developed in this section  to our motivating set of examples.

\begin{prop}
The almost-complex structure $\Om^{(\bullet,\bullet)}_q(\ccpn)$ is integrable.
\end{prop}
\demo
For  $\ol{z_{i1}} = \ol{u^i_1 S(u^1_1)} \in V\hol$, with $i = 2, \ldots, N$, we have:
\bigskip
\begin{align*}
 &  \ol{(S((z_{i1})\1))^+} \oby \ol{((z_{i1})\2)^+} \\
=  & \sum_{a = 2}^{N} \ol{(S(u^i_aS(u^b_1)))^+} \oby \ol{(u^a_1S(u^1_b))^+}  =     \sum_{a = 1}^{N} q^{2(b-1)}\ol{(u^b_1S(u^i_a))^+} \oby \ol{(u^a_1S(u^1_b))^+} \\
=  & \sum_{a = 2}^{N} \ol{(u^1_1S(u^i_a))^+} \oby \ol{(u^a_1S(u^1_1))^+}   \, +  \,  \sum_{b = 2}^{N}q^{2(b-1)} \ol{(u^b_1S(u^i_1))^+} \oby \ol{(u^1_1S(u^1_b))^+} \\
=  &  \sum_{a = 2}^{N} \ol{(u^1_1S(u^i_a))^+} \oby \ol{(u^a_1S(u^1_1))^+} = \ol{(u^1_1S(u^i_i))^+} \oby \ol{u^i_1S(u^1_1)} =  \ol{(u^1_1S(u^i_i))^+} \oby e^+_{i-1}.
\end{align*}
Thus, for any linear projection $P: \Lambda^{\oby 2} \to (V^1)^{\oby 2}$,  requirement (\ref{int:1}) of Proposition (\ref{prop-integrability}) will be satisfied.
\qed


\bigskip

Mathematical Institute of Charles University, Sokolovsk\'a 83, 186 75 Praha 8, Czech Republic

{\em e-mail}: \tt{obuachalla@karlin.mff.cuni.cz}


\begin{thebibliography}{VD} 





\bibitem{EBPS} {\sc E. Beggs, S. P. Smith}, Noncommutative complex differential geometry, {\it J. Geom. Phys.}, {\bf 72},  7-33, (2013)

\bibitem{BourAlg} {\sc N. Bourbaki}, {\it Elements of Mathematics - Algebra I},  Addison--Wesley 1974  

\bibitem{qqgauge} {\sc T. Brzezi\'nski, S. Majid}, Quantum group gauge theory on quantum spaces, {Comm. Math. Phys}, \textbf{157}, 591-638, (1993)

\bibitem{CONN} {\sc A. Connes}, {\it Noncommutative Geometry}, Academic Press, 1994

\bibitem{ConnesCuntz} {\sc A. Connes, J. Cuntz}, Quasi homomorphismes, cohomologie cyclique et positivit\'e, {\it Comm. Math. Phys.},  {\bf 114}, 515--526, (1988)

\bibitem{DDCPN} {\sc F. D'Andrea, L. D\c{a}browski}, Dirac operators on quantum projective spaces, {\it Comm. Math. Phys.}, \textbf{295}, 731-790, (2010)


\bibitem{DiengSchw} {\sc M. Dieng, A. Schwarz}, Differential and complex geometry of two-dimensional noncommutative tori, {\it Lett. Math. Phys.}, {\bf 61}, 263-270, (2002)
 

\bibitem{FGR1} {\sc J. Fr\"ochlich, O. Grandjean, A. Recknagel}, Supersymmetric quantum theory and (non-commutative) differential geometry, {\it Comm. Math. Phys.}, \textbf{193}, 527-594, (1998)


\bibitem{HK} {\sc I. Heckenberger and S. Kolb}, The locally finite part of the dual coalgebra of quantised irreducible flag manifolds, {\it Proc. Lon. Math. Soc.}, (3) \textbf{89}, 457-484, (2004)

\bibitem{HKdR} {\sc I. Heckenberger and S. Kolb}, De Rham complex for quantized irreducible flag manifolds, {\it J. Algebra}, \textbf{305}, 704-741, (2006)

\bibitem{Herm} {\sc U. Hermisson}, Derivations with quantum group action,  {\it Comm. Algebra}, {\bf 30}, 101-117, (2002)

\bibitem{HUY} {\sc D. Huybrechts}, {\it Complex geometry - An Introduction}, Springer, Universitext, 2004

\bibitem{KLVSCP1} {\sc M. Khalkhali, G. Landi, W. van Suijlekom}, Holomorphic structures on the quantum projective line, {\it Int. Math. Res. Not. IMRN,} 851-884, (2010)

\bibitem{KKCP2} {\sc M. Khalkhali, A. Moatadelro}, The homogeneous coordinate ring of the quantum projective plane, {\it J. Geom. Phys.}, {\bf 61}, 276-289, (2011)

\bibitem{KKCPN} {\sc M. Khalkhali, A. Moatadelro}, Noncommutative complex geometry of the quantum projective space, {\it J. Geom. Phys.,} {\bf 61}, 2436-2452, (2011)

\bibitem{KSLeabh} {\sc A. Klimyk, K. Schm\"udgen}, Quantum Groups and their Representations,  {\it Springer Verlag}, 1997

\bibitem{Krah} {\sc U. Kr\"ahmer}, Dirac operators on quantum flag manifolds, {\it Lett. Math. Phys.}, \textbf{67}, 49-59, (2004)

\bibitem{LR91} {\sc V. Lakshmibai, N. Reshetikhin},  Quantum deformations of flag and Schubert schemes,  {\it C.   R. Acad. Sci. Paris}, {\bf 313}, 121-126, (1991)


\bibitem{Maj1} {\sc S. Majid}, Quantum and braided group Riemannian geometry, {\it J. Geom. Phys.}, \textbf{30}, 113-146, (1999)

\bibitem{MajPrimer} {\sc S. Majid}, {\it A Quantum Groups Primer}, London Mathematical Society Lecture Note Series, {\it Cambridge University Press}, 2002

\bibitem{Maj} {\sc S. Majid}, Noncommutative Riemannian and spin geometry of the standard $q$-sphere, {\it Comm. Math. Phys.}, \textbf{256}, 255-285,
(2005)

\bibitem{Mey} {\sc U. Meyer}, Projective quantum spaces, {\it Lett. Math. Phys.}, \textbf{35}, 91-97, (1995)

\bibitem{DIAS} {\sc S. Murray, C. S\"amann}, Quantization of flag manifolds and their supersymmetric extensions, {\it Adv. Theor. Math. Phys.}, {\bf 12}, 641-710, (2008)

\bibitem{MulSch} {\sc E. F. M\"uller, H.-J. Schneider}, {Quantum homogeneous spaces with faithfully flat module structures}, {\it Israel Journal of Mathematics}, \textbf{111}, 157-190, (1999)

\bibitem{MMF1} {\sc R. \'O Buachalla}, Quantum bundle description of quantum projective spaces, {\em Comm. Math. Phys}., {\bf 316},  345-373, (2012)

\bibitem{MMF3} {\sc R. \'O Buachalla}, Noncommutative K\"ahler Structures on Quantum Homogeneous Spaces, (in preparation)

\bibitem{MMFPhD} {\sc R. \'O Buachalla}, PhD Thesis: Quantum Groups and Noncommutative Complex Geometry, London, (2013)

\bibitem{PV} {\sc G. B. Podkolzin, L.I. Vainerman}, Quantum Stiefel manifold and double cosets of quantum unitary group, {\it  Pac J. Math.}, \textbf{188}, 179-199, (1999)

\bibitem{PolishSchw} {\sc A. Polishchuk, A. Schwarz}, Categories of holomorphic vector bundles on noncommutative two-tori, {\it Comm. Math. Phys.}, \textbf{236}, 135-159, (2003)


\bibitem{Segal} {\sc G. Segal}, Equivariant K-theory, {\it Inst. Hautes \'Etudes Sci. Publ. Math.}, {\bf 34}, 129--151, (1968)

\bibitem{SV} {\sc J. T. Stafford, M. Van den Bergh}, Noncommutative curves and noncommutative surfaces,
{\it Bull. Amer. Math. Soc.,} {\bf 38}, 171-216, (2001) 

\bibitem{SoibelMirr} {\sc Y. So\v{i}bel'man}, Quantum tori, mirror symmetry and deformation theory, {\it Lett. Math. Phys.}, {\bf 56}, 99-125, (2001)

\bibitem{Soibel2} {\sc Y. S. So\v{i}bel'man}, On quantum flag manifolds, {\it Func. Ana. Appl.}, \textbf{25},
225-227, (1992)


\bibitem{Tak} {\sc M. Takeuchi}, Relative Hopf modules - equivalences and freeness conditions, {\it J. Algebra}, \text{\bf 60}, 452-471, (1979)

\bibitem{TT} {\sc E. Taft, J. Towber}, Quantum deformation of flag schemes and Grassmann schemes I. A q-deformation of the shape-algebra for
$GL(n)$, {\it J. Algebra}, {\bf 142}, 1-36, (1991)


\bibitem{Wor}
{\sc S. L. Woronowicz}, Differential calculus on compact matrix pseudogroups (quantum groups), { \it Comm. Math. Phys.}, \textbf{122},
125-170, (1989)

\end{thebibliography}
\end{document}